\documentclass[a4paper,12pt]{article}
\usepackage{xeCJK}
\usepackage{amsthm,amsmath,amssymb,mathtools}
\usepackage{mathrsfs}
\usepackage{extarrows}
\usepackage{geometry}
\usepackage{graphicx}
\usepackage{needspace}
\usepackage{color}
\usepackage{hyperref}
\hypersetup{colorlinks=true,linkcolor=blue,citecolor=blue,
filecolor=magenta,urlcolor=cyan,
pdftitle={Global existence and large-time behavior of classical axisymmetric solutions to the two-dimensional viscous Saint-Venant equations and related models},
pdfauthor={Xiangdi Huang, Muxi Lei, Weili Meng}}
\newcommand{\divg}{{\rm div}}
\newcommand{\norm}[1]{\left\Vert#1\right\Vert}

\newcommand{\dd}{\,\mathrm{d}}
\newcommand{\R}{\mathbb R}
\newtheorem{theorem}{Theorem}[section]
\newtheorem{lemma}[theorem]{Lemma}
\newtheorem{definition}[theorem]{Definition}

\newtheorem{proposition}[theorem]{Proposition}
\newtheorem{remark}[theorem]{Remark}

\numberwithin{equation}{section}
\allowdisplaybreaks
\title{Global existence and large-time behavior of classical axisymmetric
solutions to the two-dimensional viscous Saint-Venant equations
and related models}
\author{Xiangdi Huang$^{a}$, Muxi Lei$^{a}$, Weili Meng$^{a}$%
\thanks{Corresponding author: Weili Meng. E-mail addresses:
\href{mailto:xdhuang@amss.ac.cn}{xdhuang@amss.ac.cn} (X. Huang);
\href{mailto:leimuxi25@mails.ucas.ac.cn}{leimuxi25@mails.ucas.ac.cn} (M. Lei);
\href{mailto:mengweili@amss.ac.cn}{mengweili@amss.ac.cn} (W. Meng).}\\[0.4em]
{\small $^{a}$ Institute of Mathematics, Academy of Mathematics and Systems Science,}\\
{\small Chinese Academy of Sciences, Beijing 100190, China}}
\date{}
\begin{document}
\maketitle
\begin{abstract}
In this paper, we study a two-dimensional compressible Navier--Stokes system
with density-dependent viscosity in a disk, subject to the no-slip boundary
condition. We consider both the symmetric viscous stress and the stress
written in terms of the full velocity gradient. In each case the
coefficients satisfy the BD relation $\mu(\rho)=\rho^\alpha$ and
$\lambda(\rho)=2(\alpha-1)\rho^\alpha$. For every $1/2<\alpha\le1$ and
$\gamma>1$, we establish the global existence and uniqueness of classical
solutions for arbitrarily large axisymmetric initial data with strictly
positive density and possibly nonzero tangential velocity.
Moreover, the density has upper and positive
lower bounds independent of time, and the solution converges to the
constant equilibrium. The two stresses give the same radial momentum
equation, but require different estimates of the tangential velocity.
The endpoint $\alpha=1$, $\gamma=2$ includes the viscous shallow-water
equations.


\end{abstract}
\tableofcontents

\section{Introduction}

In this paper, we consider the following compressible system in
$\Omega=B_R\subset\R^2$:
\begin{equation}\label{eq:model}
\left\{
\begin{aligned}
&\rho_t+\divg(\rho u)=0,\\
&\rho(u_t+u\cdot\nabla u)+\nabla P(\rho)
=\divg(2\mu(\rho)\mathcal T u)+\nabla(\lambda(\rho)\divg u).
\end{aligned}
\right.
\end{equation}
We consider the two choices
\begin{equation}\label{eq:stress}
\mathcal T u=\mathbb Du+\sigma\mathbb Au,\qquad
\mathbb Du=\frac{\nabla u+(\nabla u)^\top}{2},\qquad
\mathbb Au=\frac{\nabla u-(\nabla u)^\top}{2},
\qquad \sigma\in\{0,1\}.
\end{equation}
Thus $\sigma=0$ gives the symmetric stress and $\sigma=1$ gives
$\mathcal T u=\nabla u$, which we call the gradient stress.
Here $\rho$ and $u$ denote the density and velocity, respectively, and
$P(\rho)=\rho^\gamma$ with $\gamma>1$. The viscosity coefficients are
\begin{equation}\label{eq:viscosity}
\mu(\rho)=\rho^\alpha,\qquad
\lambda(\rho)=2(\alpha-1)\rho^\alpha,\qquad
\frac12<\alpha\le1.
\end{equation}
They satisfy the algebraic relation
$\lambda(\rho)=2(\rho\mu'(\rho)-\mu(\rho))$.
The initial and boundary conditions are
\begin{equation}\label{eq:data}
(\rho,u)(x,0)=(\rho_0,u_0)(x),\qquad
u|_{\partial\Omega}=0.
\end{equation}

Controlling the density is a central difficulty in the large-data theory of
compressible Navier--Stokes equations. Density concentration can obstruct
continuation even for nondegenerate viscosities. When the viscosities
depend on the density and vanish at vacuum, a positive lower bound is also
needed to retain the ellipticity of the momentum equation. We review the
results most closely related to these two issues.

For constant viscosity coefficients, the classical local theory away from
vacuum goes back to Nash's existence result \cite{Nash1962} and Serrin's
uniqueness theory \cite{Serrin1959}.
Cho and Kim \cite{CK2006} developed a local classical theory allowing
vacuum under compatibility conditions. Matsumura and Nishida
\cite{MN1980} established
global existence for small smooth perturbations of a positive equilibrium,
and Danchin \cite{Danchin2000} developed the corresponding theory in critical
spaces for the barotropic equations. Hoff \cite{Hoff1995} constructed global
solutions for small perturbations that may have discontinuities. For the
three-dimensional isentropic system, Huang, Li and Xin \cite{HLX2012}
obtained global classical solutions with small initial energy, allowing
large oscillations and vacuum. Li and Xin \cite{LiXin2019} established
global well-posedness and large-time behavior for the corresponding
two-dimensional Cauchy problem with small initial energy and vacuum at infinity.

For general large data, Lions \cite{Lions1998} constructed global
finite-energy weak solutions, in particular for $\gamma\ge9/5$ in three
dimensions. Feireisl, Novotn\'y and Petzeltov\'a \cite{FNP2001} relaxed
the three-dimensional condition to $\gamma>3/2$. Under spherical symmetry,
Jiang and Zhang \cite{JZ2001} obtained global weak solutions for
$\gamma>1$; see also their work on three-dimensional axisymmetric
flows \cite{JZ2003}. For smooth spherically symmetric flows in solid
balls, Huang and Matsumura \cite{HuangMatsumura2015} characterized
finite-time breakdown in terms of density concentration or vacuum near
the center. Merle, Rapha\"el, Rodnianski and Szeftel \cite{MRRS2022}
constructed finite-time implosions for three-dimensional compressible
flows in a suitable range of pressure laws.

A different approach to large initial data is to keep the shear
viscosity constant and take the bulk viscosity proportional to
$\rho^\beta$. Vaigant and Kazhikhov \cite{VK1995} obtained global
solutions in two dimensions for $\beta>3$ and positive initial density.
On the torus, Jiu, Wang and Xin \cite{JWX2014} treated initial vacuum
for $\beta>3$, and Huang and Li \cite{HL2016} lowered the condition to
$\beta>4/3$. The whole-plane problem with far-field vacuum was treated
in weighted spaces by Jiu, Wang and Xin \cite{JWX2013} and
by Huang and Li \cite{HL2022}, who established global classical solutions
for $\beta>4/3$ under suitable decay and regularity conditions.
Fan, Jiang, Li and L\"u \cite{FJLL2026Cauchy} proved global existence
and uniqueness of strong solutions for large data allowing vacuum in both
settings for every $\beta>1$. Their result also covers near-linear laws
such as $\lambda(\rho)=\rho[\log(e+\rho)]^s$ with $s>1/2$. They obtained
exponential convergence to equilibrium on the torus and algebraic decay
of the density toward vacuum in the whole plane.

In bounded domains, the boundary condition affects the estimates for
the effective viscous flux. In simply connected domains with Navier-slip boundary conditions,
Fan, Li and Li \cite{FLL2022} established global solvability for large
data allowing vacuum when $\beta>4/3$. Fan, Li and Wang \cite{FLW2023}
studied convergence to equilibrium in this setting, including the case
of large potential forces. In smooth bounded domains, possibly with holes,
Fan et al. \cite{FJLL2026Slip} obtained global strong solutions and
exponential convergence to equilibrium for every $\beta>1$ when the
curvature-corrected slip coefficient is nonnegative and not identically
zero. Under the no-slip condition in a disk,
Huang, Su, Yan and Yu \cite{HSYY2024} proved global existence of
radial strong solutions without swirl for $\beta>1$. Huang and Meng
\cite{HM2025} treated axisymmetric flows under the same boundary condition, obtaining
global strong solutions and their large-time behavior for $\beta\ge1$
with initial vacuum allowed. They also treated $0<\beta<1$ under an
additional upper bound on the initial density.

For the shallow-water viscosity law, Bresch, Desjardins and Lin
\cite{BDL2003} derived an entropy estimate that controls a density gradient.
Bresch and Desjardins \cite{BD2004} identified the general algebraic relation
underlying this estimate, $\lambda(\rho)=2(\rho\mu'(\rho)-\mu(\rho))$,
and later used the estimate to construct approximate solutions \cite{BD2006}.
Bresch, Desjardins and G\'erard-Varet
\cite{BDGV2007} used the entropy structure to construct weak solutions in
bounded domains for a model with drag and suitable boundary conditions.
Mellet and Vasseur \cite{MV2007} established compactness of weak solutions using improved
velocity integrability. Guo, Jiu and Xin \cite{GJX2008} constructed
global spherically symmetric weak solutions. For flows without a symmetry
assumption, Li and Xin \cite{LiXinBD2015} proved global weak existence
for classes of degenerate viscosities, and Vasseur and Yu \cite{VY2016}
treated the three-dimensional case with linear shear viscosity and zero
bulk viscosity.
Bresch, Vasseur and Yu \cite{BVY2022} obtained entropy-weak solutions
for nonlinear density-dependent viscosities with the physical symmetric
stress.

For global regular solutions in two or three spatial dimensions, the BD
entropy estimate alone does not yield pointwise upper and positive lower
bounds on the density. Under spherical symmetry,
Zhang \cite{Zhang2025} obtained global strong solutions with large
positive-density data in balls. Huang, Meng and Zhang
\cite{HMZ2025} extended the parameter range for large-data classical
solutions, including the two-dimensional viscous shallow-water endpoint,
and obtained time-uniform density bounds and large-time behavior in
specified parameter ranges. For the shallow-water viscosity law, Chen,
Zhang and Zhu \cite{CZZ2025} obtained global regular solutions of the
isentropic Cauchy problem in two and three dimensions, allowing far-field
vacuum. Gu and Huang \cite{GH2025} treated the same viscosity law with
transported entropy and proved global classical solvability for large
spherically symmetric data away from vacuum in a disk for $\gamma>1$ and
in a three-dimensional ball for $1<\gamma<3$. For sublinear viscosity powers, Lei \cite{LeiBD2026} established
global classical solutions in the whole space with positive far-field
density. Chen, Zhang and Zhu \cite{CZZ2026} subsequently treated
sublinear powers with far-field vacuum. For a range of superlinear viscosity powers, Huang
and B.~Li \cite{HLBD2026} identified an initial boundary sign condition
on the effective velocity for global existence in a three-dimensional
solid ball and constructed examples of breakdown by boundary vacuum.
For classical solutions in balls, Lei and Liang \cite{LL2026} extended
the large-time convergence result of Huang, Meng and Zhang \cite{HMZ2025}
to a wider range of viscosity and pressure exponents.

We study axisymmetric flows in a disk, as in \cite{HM2025}, with the
density-dependent shear viscosity of \cite{HMZ2025}. In the purely radial
setting of \cite{HMZ2025}, the antisymmetric part of the velocity gradient
vanishes. With swirl, it contributes to the vector effective-velocity
equation, and the radial equation acquires a centrifugal forcing term.
Controlling this term without assuming a positive lower bound for the
density is the main step in the argument below. 

The stresses in \eqref{eq:stress} agree for radial flows without swirl.
With swirl, their difference is tangential after taking the divergence
in \eqref{eq:model}. The radial effective-velocity equation is therefore
the same for both stresses. The angular equations, however, are
different, and we treat their estimates separately.

\subsection{Notations and main result}

For $k\ge0$ and $1\le p\le\infty$, we use the standard spaces
$L^p=L^p(\Omega)$, $W^{k,p}=W^{k,p}(\Omega)$ and
$H^k=W^{k,2}(\Omega)$. We write $\norm f_p$ for the $L^p$ norm.
The material derivative is $\dot f=f_t+u\cdot\nabla f$, and
$\bar f=|\Omega|^{-1}\int  f\dd x$.
For a scalar $a$, let $a_+=\max\{a,0\}$ and $a_-=\max\{-a,0\}$.

We call $(\rho,u)$ axisymmetric if, for every planar rotation $Q$,
$\rho(Qx,t)=\rho(x,t)$ and $u(Qx,t)=Q u(x,t)$. Equivalently,
\begin{equation}\label{eq:symmetry}
\rho=\rho(r,t),\qquad
u=U(r,t)e_r+W(r,t)e_\theta,\qquad r=|x|,
\end{equation}
where $e_r=(x_1,x_2)^\top/r$ and
$e_\theta=(x_2,-x_1)^\top/r$. Here $U$ is the radial velocity and $W$
is the tangential velocity (swirl). The purely radial case is $W\equiv0$.
All spatial integrals without specified limits
are over $\Omega$ with respect to the two-dimensional measure $\dd x$.
In particular, for a radial integrand,
\begin{equation}\label{eq:measure}
\int  f(r)\dd x=2\pi\int_0^Rf(r)r\dd r.
\end{equation}
The effective velocity is defined by
\begin{equation}\label{eq:effective}
v=u+2\alpha\rho^{\alpha-2}\nabla\rho
=Ve_r+We_\theta,\qquad
V=U+2\alpha\rho^{\alpha-2}\rho_r.
\end{equation}

\begin{definition}[Global classical solutions]\label{def:classical}
An axisymmetric pair $(\rho,u)$ with $\rho>0$
is called a global classical solution if it satisfies
\eqref{eq:model}--\eqref{eq:data} and, for every $0<\tau<T<\infty$,
\[
\left\{
\begin{aligned}
&\rho\in C([0,T];H^3),\qquad \rho_t\in C([0,T];H^2),\\
&\rho_{tt}\in L^\infty(0,T;L^2)\cap L^2(0,T;H^1),\\
&u\in C([0,T];H^3\cap H_0^1)\cap L^2(0,T;H^4),\\
&u_t\in L^\infty(0,T;H_0^1)\cap L^2(0,T;H^2),\qquad
u_{tt}\in L^2(0,T;L^2),\\
&u\in L^\infty(\tau,T;H^4),\qquad
u_t\in L^\infty(\tau,T;H^2),\\
&u_{tt}\in L^\infty(\tau,T;L^2)\cap L^2(\tau,T;H_0^1).
\end{aligned}
\right.
\]
\end{definition}

\begin{theorem}\label{thm:main}
Assume that $1/2<\alpha\le1$ and $\gamma>1$, and fix either choice of
$\mathcal T$ in \eqref{eq:stress}. Let the initial data be axisymmetric
in the sense of \eqref{eq:symmetry} and satisfy
\begin{equation}\label{eq:initial}
0<c_0\le\rho_0\le C_0,\qquad
(\rho_0,u_0)\in H^3,\qquad u_0|_{\partial\Omega}=0,
\end{equation}
and the compatibility condition
\begin{equation}\label{eq:compatibility}
\rho_0^{-1}\left[
\divg(2\rho_0^\alpha\mathcal T u_0)
+\nabla\bigl(2(\alpha-1)\rho_0^\alpha\divg u_0\bigr)
-\nabla\rho_0^\gamma\right]
-u_0\cdot\nabla u_0\in H_0^1.
\end{equation}
Then the initial-boundary-value problem
\eqref{eq:model}--\eqref{eq:data} admits a unique global classical solution.
Moreover, there is a constant $C>0$, independent of time, such that
\begin{equation}\label{eq:main-density}
C^{-1}\le\rho(x,t)\le C\qquad
\text{for }(x,t)\in\overline\Omega\times[0,\infty),
\end{equation}
Furthermore, we have
\begin{equation}\label{eq:main-asymptotic}
\norm{\rho(t)-\bar\rho}_\infty+
\norm{\nabla\rho(t)}_2+
\norm{\nabla u(t)}_2+
\norm{\nabla^2u(t)}_2\longrightarrow0
\quad\text{as }t\to\infty.
\end{equation}
The constants may depend on $\alpha,\gamma,R$ and the initial data.
\end{theorem}

\begin{remark}\label{rmk:shallow-water}
The choice $\alpha=1$ and $\gamma=2$ with $\sigma=0$ corresponds to
the viscous shallow-water equations with symmetric stress. Thus Theorem
\ref{thm:main} includes global existence and large-time behavior for
axisymmetric classical solutions of this model with positive density.
The corresponding gradient-stress system is also covered.
\end{remark}

\begin{remark}\label{rmk:stress}
The compatibility condition \eqref{eq:compatibility} uses the chosen
stress in \eqref{eq:stress}. The two conclusions of Theorem \ref{thm:main}
concern separate initial-boundary value problems. In the presence of
swirl, one cannot replace $\mathbb Du$ by $\nabla u$ in the momentum
equation.
\end{remark}

\begin{remark}\label{rmk:comparison}
Under spherical symmetry without swirl, Huang--Meng--Zhang
\cite{HMZ2025} established the global existence of classical solutions
for $1/2<\alpha<1$ and $\gamma>1$, and for $\alpha=1$ and
$\gamma\ge3/2$. They also obtained time-uniform upper and positive lower 
bounds for the density and the large-time behavior of the solutions
for $0.54<\alpha<1$ and $\gamma>1$, and for $\alpha=1$ and
$\gamma\ge3/2$. Theorem \ref{thm:main} extends their
result to axisymmetric flows.
\end{remark}

\begin{remark}\label{rmk:effective-velocity}
In the BD approach used in Huang--Meng--Zhang \cite{HMZ2025}, the
effective velocity plays a key role in estimating the density.
For the system \eqref{eq:model}, the effective velocity
$v=u+2\alpha\rho^{\alpha-2}\nabla\rho$ satisfies
\begin{equation}\label{eq:effective-system}
\rho v_t+\rho u\cdot\nabla v+\nabla P
=2(1+\sigma)\divg(\rho^\alpha\mathbb Au).
\end{equation}
For radial flows without swirl, $\mathbb Au=0$, and the right-hand
side vanishes. This allows one to derive the density-weighted
$L^p$ estimates of $v$ for suitable $p>2$, which are essential for
the density bounds. When swirl is present, $\mathbb Au$ need not
vanish. Obtaining higher integrability of the effective velocity
and time-uniform upper and positive lower bounds for the density
is therefore the main difficulty in this paper.
\end{remark}

\subsection{Main ideas of the proof}

In the two-dimensional case studied by Huang--Meng--Zhang \cite{HMZ2025},
the density bounds rely on estimates of
$\nabla\rho^{\alpha-1+1/m}$ in $L^\infty(0,T;L^m)$ for suitable
$m>2$. Indeed, the definition of the effective velocity gives
\[
 \partial_r\rho^{\alpha-1+1/m}
 =\frac{\alpha-1+1/m}{2\alpha}\rho^{1/m}(V-U).
\]
Thus it suffices to estimate $\rho^{1/m}u$ and $\rho^{1/m}v$ in
$L^\infty(0,T;L^m)$. Without swirl, the
antisymmetric stress in \eqref{eq:effective-system} vanishes, and one
can test the resulting transport equation by powers of the effective
velocity. With swirl, the same test encounters derivatives of the
antisymmetric part of $\nabla u$. Our first observation is that
$\nabla\rho$ is radial, so only the radial component $V$ is needed.
The divergence of $\rho^\alpha\mathbb Au$ is tangential under
\eqref{eq:symmetry}, whereas the radial component of $u\cdot\nabla v$
contains $-W^2/r$. The radial equation \eqref{eq:radial-effective}
therefore reads
\[
 \rho(V_t+UV_r)+\frac{\gamma}{2\alpha}\rho^{\gamma-\alpha+1}V
 =\frac{\gamma}{2\alpha}\rho^{\gamma-\alpha+1}U
  +\frac{\rho W^2}{r}.
\]
This reduces the problem to estimating a scalar centrifugal term.
Its sign is favorable for $V_-$, but a bound for $\rho^{1/m}V$
in $L^\infty(0,T;L^m)$ requires a space-time estimate of $W^2/r$
near the center.

We first discuss the gradient stress. When $1/2<\alpha<1$,
Proposition \ref{oa:velocity}
gives $\sup_{0\le t\le T}\int\rho|u|^m\,dx\le C_m$ for
$2<m<N_\alpha$, where $N_\alpha$ is given by \eqref{oa:N}. Using the sign of the centrifugal term,
we estimate $V_-$ at the order $m=1/(1-\alpha)$, which allows us
to control the negative part of $(\log\rho)_r$ in $L^m$. Integrating in $r$ gives the one-sided comparison
\[
 \rho(r,t)\le C\rho(s,t),\qquad 0\le r\le s\le R.
\]
Together with the bound for $\rho(R,t)$, this gives an upper bound
for the density. The same comparison with $r=0$ reduces the lower
bound to controlling the density at the center.

To estimate $W^2/r$, we test \eqref{eq:angular} by
$|W|^{q-2}Wr^{-k}$. The resulting dissipation is
\[
 2\int\rho^\alpha|W|^{q-2}r^{-k}
 \left((q-1)W_r^2-kW_r\frac Wr+\frac{W^2}{r^2}\right)\,dx.
\]
The last square gives the desired weight $r^{-M}$ when
$q=2M$ and $k=M-2$. The quadratic form is pointwise positive definite
for $k^2<4(q-1)$, which for this choice gives $2<M<6+2\sqrt7$.
Since $\rho$ is bounded above, this gives a bound for
$\int_0^T\int\rho(W^2/r)^M\,dx\,dt$. The $V$ equation needs
a bound for $\int_0^T(\int\rho(W^2/r)^m\,dx)^{1/m}\,dt$.
To obtain it, we take $M>m$ and interpolate with
$\int\rho W^2/r\,dx$, which decays exponentially by the energy decay
and a low-order weighted angular estimate. This gives a uniform
$L^m$ bound for $\rho^{1/m}V$ for $2<m<\min\{N_\alpha,6+2\sqrt7\}$.

To obtain a lower bound by comparing density values, we choose
$m>1/(1-\alpha)$ so that $\rho^{\alpha-1+1/m}$ is a negative
power of $\rho$. This choice is compatible
with $m<\min\{N_\alpha,6+2\sqrt7\}$ when
\[
 \frac12<\alpha<1-\frac1{6+2\sqrt7}=\frac{1+\sqrt7}{4},
\]
since $N_\alpha>1/(1-\alpha)$ throughout $1/2<\alpha<1$.
In this range, the oscillation bound for the negative power gives a
lower bound for the central density by comparison with a point where
$\rho=\bar\rho$. For $(1+\sqrt7)/4\le\alpha<1$, however,
$1/(1-\alpha)\ge6+2\sqrt7$, so no $m$ in the range given by the
direct angular estimate satisfies $m>1/(1-\alpha)$.

To reach the required range of $m$, we freeze the viscosity coefficient at its
central value $\rho(0,t)^\alpha$. This recovers positivity of the
integrated angular dissipation near the center even when the
quadratic form is not pointwise positive. For $1/(1-\alpha)<m<N_\alpha$,
the resulting bound for $\rho^{1/m}V$ in $L^\infty(0,T;L^m)$ depends on
a power of $H=\max\{1,\sup_{0\le t\le T}\rho(0,t)^{-1}\}$.
A finite iteration reduces this power to some
$\varepsilon<1-\alpha-1/m$. The estimate for $V-U$ then bounds
the oscillation of $\rho^{\alpha-1+1/m}$ by $CH^\varepsilon$.
Comparing this negative power at the center with its value at a
point where $\rho=\bar\rho$ gives
\[
 H^{1-\alpha-1/m}\le C(1+H^\varepsilon).
\]
The one-sided comparison then gives
a uniform positive lower bound throughout the disk for every
$1/2<\alpha<1$.

At $\alpha=1$, we instead control
$(\log\rho)_r=(V-U)/2$ by bounds for $U$ and $V$ on finite time
intervals. After obtaining the density upper bound, we establish
lower bounds for $U$ and $V$ by a coupled comparison. These lower
bounds allow us to bound $|W|$ and $U_+$ before estimating $V_+$. The remaining issue is the factor $1/r$ in the centrifugal term.
A comparison for $W/\sqrt r$ near the center yields
\[
 \left\|\frac{W(t)^2}{r}\right\|_{L^\infty}
 \le C(T)\left(1+\sup_{0\le s\le t}\|V_+(s)\|_{L^\infty}\right),
 \qquad 0\le t\le T.
\]
The key is that the right-hand side is linear in the running
supremum of $V_+$. Thus Gronwall's inequality closes the estimate
for $V_+$, and the bound for $(\log\rho)_r$ gives finite-time
positivity and global continuation.

To obtain a lower bound independent of time at $\alpha=1$, we combine
weighted angular estimates and energy decay to control the centrifugal
term in the estimate for $\int\rho|V|^3\,dx$. The resulting uniform
bound gives a uniform $L^3$ bound for $\nabla\rho$. Interpolating
with the $L^2$ decay of $\rho-\bar\rho$ gives uniform convergence to
the positive mean density. The large-time lower bound and the lower
bound on the remaining finite interval then give a uniform positive
lower bound for all time.

For the symmetric stress, the angular shear $W_r-W/r=r(W/r)_r$
measures the radial variation of the angular velocity $W/r$.
To control the centrifugal term $W^2/r$, we rewrite the angular
equation using $2(\rho^\alpha)_r=\rho(V-U)$. We first bound $U$ and
$V$ from below by comparison. The lower bound on $V$, together with
the enclosed mass estimate and the BD bound, gives $\rho\le C$.
The maximum principle applied to
$W/\sqrt r$ then bounds $W^2/r$ on each finite time interval.
The radial equations now bound $U$ and $V$ on that interval, and
integration of \eqref{eq:effective} gives the positive lower density
bound needed for continuation.

To make the lower bound uniform in time, we combine weighted estimates
from the two forms of the angular equation. The resulting dissipation
controls powers of $W/r$ near the center. On each fixed outer annulus,
energy decay and the BD estimate give uniform convergence of $\rho$ to
its positive mean. Using this annular lower bound and the integrability
of the relative energy, we obtain
$\int_0^\infty \|\rho^{1/m}W^2/r\|_m\dd t<\infty$ for some $m>2$.
The radial equations then give uniform density-weighted $L^m$ bounds
for $U$ and $V$, so \eqref{eq:effective} bounds $\nabla\rho$ in $L^m$.
Interpolation with $\|\rho-\bar\rho\|_2\to0$ yields uniform convergence
throughout the disk and hence a time-independent positive lower bound.

Once these bounds are available, the higher-order estimates and the
large-time behavior follow by the arguments of Huang--Meng--Zhang
\cite[Sections 3.6--3.8]{HMZ2025}.

\section{Preliminary estimates}

We first state the local existence result for positive initial density.
It can be established by adapting the arguments in
\cite{CK2006,HMZ2025} to the system \eqref{eq:model}--\eqref{eq:data}.

\begin{lemma}[Local existence of classical solutions]\label{lem:local}
Assume that $1/2<\alpha\le1$ and $\gamma>1$, and that the axisymmetric
initial data $(\rho_0,u_0)$ satisfy \eqref{eq:initial} and
\eqref{eq:compatibility}. Then there exists a time $T_0>0$ such that
\eqref{eq:model}--\eqref{eq:data} admits a unique axisymmetric classical
solution $(\rho,u)$ on $\Omega\times[0,T_0]$, with
\[
\frac{c_0}{2}\le\rho(x,t)\le2C_0
\qquad\text{for }(x,t)\in\overline\Omega\times[0,T_0].
\]
Moreover, for every $0<\tau<T_0$, the solution satisfies
\[
\left\{
\begin{aligned}
&\rho\in C([0,T_0];H^3),\qquad \rho_t\in C([0,T_0];H^2),\\
&\rho_{tt}\in L^\infty(0,T_0;L^2)\cap L^2(0,T_0;H^1),\\
&u\in C([0,T_0];H^3\cap H_0^1)\cap L^2(0,T_0;H^4),\\
&u_t\in L^\infty(0,T_0;H_0^1)\cap L^2(0,T_0;H^2),\qquad
u_{tt}\in L^2(0,T_0;L^2),\\
&u\in L^\infty(\tau,T_0;H^4),\qquad
u_t\in L^\infty(\tau,T_0;H^2),\\
&u_{tt}\in L^\infty(\tau,T_0;L^2)\cap L^2(\tau,T_0;H_0^1).
\end{aligned}
\right.
\]
\end{lemma}

By Lemma \ref{lem:local}, there is a unique positive local classical
solution. Let $T_*$ be its maximal existence time and fix $0<T<T_*$.
We shall first derive estimates which do not depend on the minimum of
$\rho$ on $\overline\Omega\times[0,T]$, and then use these estimates
to establish upper and positive lower bounds for the density.
Throughout the proof, $C$ denotes a positive constant depending on
$\alpha,\gamma,R$ and the initial data, but independent of $T$.
We write $C(T)$ when dependence on $T$ is allowed. Constants with a fixed
exponent as a subscript may also depend on that exponent.

To express the velocity estimates below, we use
\begin{equation}\label{eq:Dp}
D_p(t)=\int\rho^\alpha|u|^p|\nabla u|^2\dd x,\qquad p\ge0.
\end{equation}
For a field of the form \eqref{eq:symmetry}, direct differentiation gives
\begin{equation}\label{eq:polar-gradient}
|\nabla u|^2=U_r^2+W_r^2+\frac{U^2+W^2}{r^2},
\qquad \divg u=U_r+\frac Ur.
\end{equation}
The regularity at the center and the no-slip condition imply
$U(0,t)=W(0,t)=U(R,t)=W(R,t)=0$. The radial effective velocity also
satisfies $V(0,t)=0$.

We first derive the basic energy and BD entropy estimates.
\begin{lemma}[Energy and BD entropy]\label{lem:basic}
There exists a constant $C$, independent of $T$, such that
\begin{equation}\label{eq:basic}
\sup_{0\le t\le T}\int\left(\rho|u|^2+\rho^\gamma\right)\dd x
+\int_0^TD_0(t)\dd t\le C
\end{equation}
and
\begin{equation}\label{eq:bd}
\sup_{0\le t\le T}\int\rho|v|^2\dd x+
\int_0^T\int\left(
\rho^{\gamma+\alpha-3}|\nabla\rho|^2+
\rho^\alpha|\mathbb Au|^2\right)\dd x\dd t\le C,
\end{equation}
In particular,
\begin{equation}\label{eq:sphere}
\sup_{0\le t\le T}\norm{\nabla\rho^{\alpha-1/2}}_2\le C,\qquad
\sup_{0\le t\le T}\norm{\rho^a r^b}_\infty\le C_{a,b}
\quad\text{for every }a,b>0.
\end{equation}
Moreover, $\rho(R,t)\le C$ for $0\le t\le T$.
\end{lemma}
\begin{proof}
Multiplying the second equation in \eqref{eq:model} by $u$ and using
the first equation and the boundary condition \eqref{eq:data}, we obtain
\begin{equation}\label{eq:energy-identity}
\frac{\dd}{\dd t}\int\left(\frac12\rho|u|^2+
\frac{\rho^\gamma}{\gamma-1}\right)\dd x
+2\int\rho^\alpha\left(|\mathbb Du|^2+\sigma|\mathbb Au|^2
 +(\alpha-1)(\divg u)^2\right)\dd x=0.
\end{equation}
Since $(\divg u)^2\le2|\mathbb Du|^2$ and $\alpha\le1$, the
dissipation in \eqref{eq:energy-identity} controls
$\int\rho^\alpha|\mathbb Du|^2\,dx$ for either stress, and also
$\int\rho^\alpha|\mathbb Au|^2\,dx$ when $\sigma=1$.

Next, testing \eqref{eq:effective-system} by $v$, and using
$\mathbb Av=\mathbb Au$ and the first equation in \eqref{eq:model},
we obtain
\begin{align*}
&\frac{\dd}{\dd t}\int\left(\frac12\rho|v|^2+
\frac{\rho^\gamma}{\gamma-1}\right)\dd x\\
&\qquad+2\alpha\gamma\int\rho^{\gamma+\alpha-3}|\nabla\rho|^2\dd x
+2(1+\sigma)\int\rho^\alpha|\mathbb Au|^2\dd x=0.
\end{align*}
The boundary term vanishes because $v$ is normal to $\partial\Omega$
and $(\mathbb Au)n\cdot n=0$. Integrating this identity proves
\eqref{eq:bd}. Combining its dissipation with
\eqref{eq:energy-identity} and using
$|\nabla u|^2=|\mathbb Du|^2+|\mathbb Au|^2$ gives
\eqref{eq:basic} for both stresses.
By \eqref{eq:effective}, \eqref{eq:basic} and \eqref{eq:bd},
\[
\int|\nabla\rho^{\alpha-1/2}|^2\dd x
\le C_\alpha\int\rho(|v|^2+|u|^2)\dd x\le C.
\]

Mass conservation and $0<2\alpha-1\le1$ also give
$\norm{\rho^{\alpha-1/2}}_2\le C$. Thus $\rho^{\alpha-1/2}$ is
bounded in $H^1(\Omega)$, and the one-dimensional trace inequality
on $[R/2,R]$ yields $\rho(R,t)\le C$.
For $0<r<R$, radial Cauchy--Schwarz gives
\begin{align*}
\rho(r,t)^{\alpha-1/2}
&\le\rho(R,t)^{\alpha-1/2}
+\left(\int_r^R|\partial_s\rho^{\alpha-1/2}|^2s\dd s\right)^{1/2}
\left(\log\frac Rr\right)^{1/2}\\
&\le C\left(1+\left|\log\frac rR\right|^{1/2}\right).
\end{align*}
Every positive power of this logarithm is bounded after multiplication
by $r^b$, for every $b>0$. This proves \eqref{eq:sphere}.
\end{proof}

For both stresses, the radial component of
\eqref{eq:effective-system} gives
\begin{equation}\label{eq:radial-effective}
\rho(V_t+UV_r)+\frac{\gamma}{2\alpha}
\rho^{\gamma-\alpha+1}V
=\frac{\gamma}{2\alpha}\rho^{\gamma-\alpha+1}U+\frac{\rho W^2}{r},
\end{equation}
For the gradient stress, the tangential component of the second equation
in \eqref{eq:model} gives
\begin{equation}\label{eq:angular}
\rho(W_t+UW_r)+\frac{\rho UW}{r}
=\frac2r(r\rho^\alpha W_r)_r-\frac{2\rho^\alpha W}{r^2}.
\end{equation}

We define the relative energy by
\begin{equation}\label{eq:relative-energy}
\mathcal E(t)=\int\left\{\frac12\rho|u|^2+
\frac{\rho^\gamma-\bar\rho^\gamma
-\gamma\bar\rho^{\gamma-1}(\rho-\bar\rho)}{\gamma-1}\right\}\dd x.
\end{equation}

For the gradient stress, we next prove exponential decay of the energy
under an additional upper bound for the density. We use this estimate
later to obtain a uniform positive lower bound for $\rho$.
\begin{lemma}[Exponential decay for the gradient stress]\label{lem:decay}
Let $\sigma=1$ and assume in addition that $\rho\le M$ on $\overline\Omega\times[0,T]$,
where $M$ is independent of $T$. Then there exist constants
$C,\lambda>0$, depending on $M,\alpha,\gamma,R$ and the initial data,
but independent of $T$ and of the positive lower bound of $\rho$,
such that, for every $0\le t\le T$,
\begin{equation}\label{eq:relative-decay}
\int\rho|u(t)|^2\dd x+\norm{\rho(t)-\bar\rho}_2^2
+\int_t^TD_0(s)\dd s\le Ce^{-\lambda t}.
\end{equation}
\end{lemma}
\begin{proof}
Since $\bar\rho>0$ and $0\le\rho\le M$, the strict convexity of
$s\mapsto s^\gamma$, together with its Taylor expansion at
$s=\bar\rho$, implies
\begin{equation}\label{eq:relative-equivalence}
\begin{split}
C^{-1}\left(\int\rho|u|^2\dd x+\norm{\rho-\bar\rho}_2^2\right)
&\le\mathcal E(t)\le C\left(\int\rho|u|^2\dd x+\norm{\rho-\bar\rho}_2^2\right).
\end{split}
\end{equation}
Indeed, the ratio of the pressure potential to $(\rho-\bar\rho)^2$
extends to a continuous strictly positive function on $[0,M]$.
By mass conservation and \eqref{eq:energy-identity},
\begin{equation}\label{eq:relative-energy-inequality}
\mathcal E'(t)+2(2\alpha-1)D_0(t)\le0.
\end{equation}

We next construct the radial test field
\[
\mathcal B(x,t)=\frac{e_r}{r}\int_0^r
s\bigl(\rho(s,t)-\bar\rho\bigr)\dd s.
\]
Its definition and mass conservation give
\[
\divg\mathcal B=\rho-\bar\rho,\qquad
\mathcal B|_{\partial\Omega}=0,\qquad
\mathcal B(0,t)=0.
\]
Writing the first equation in \eqref{eq:model} as
$\rho_t+r^{-1}(r\rho U)_r=0$, we have
\begin{align*}
\mathcal B_t
&=\frac{e_r}{r}\int_0^rs\rho_t(s,t)\dd s
=-\frac{e_r}{r}\int_0^r\partial_s(s\rho U)\dd s
=-\rho Ue_r.
\end{align*}
Moreover, since $0<\bar\rho\le M$,
\[
\left|\frac{\mathcal B\cdot e_r}{r}\right|
\le\frac1{r^2}\int_0^rs|\rho(s,t)-\bar\rho|\dd s\le M,
\qquad
\partial_r(\mathcal B\cdot e_r)
=\rho-\bar\rho-\frac{\mathcal B\cdot e_r}{r}.
\]
Thus $\norm{\nabla\mathcal B}_\infty\le C(M)$.
Since $\mathcal B\cdot e_r$ vanishes at both $r=0$ and $r=R$,
integration by parts also yields
\begin{align}
\norm{\nabla\mathcal B}_2^2
&=2\pi\int_0^R
\left\{|\partial_r(\mathcal B\cdot e_r)|^2
+\frac{|\mathcal B\cdot e_r|^2}{r^2}\right\}r\dd r\notag\\
&=2\pi\int_0^R
\left|\partial_r(\mathcal B\cdot e_r)
+\frac{\mathcal B\cdot e_r}{r}\right|^2r\dd r
=\norm{\rho-\bar\rho}_2^2.                 \label{eq:testfield-gradient}
\end{align}
Here the difference between the first two integrals is the boundary
term $2\pi[|\mathcal B\cdot e_r|^2]_0^R=0$.

Set $I(t)=\int\rho u\cdot\mathcal B\dd x$.
Using both equations in \eqref{eq:model}, and integrating by parts
with $\mathcal B|_{\partial\Omega}=0$, we obtain
\begin{align}
I'(t)
={}&\int(\rho u)_t\cdot\mathcal B\dd x
+\int\rho u\cdot\mathcal B_t\dd x\notag\\
={}&\int\rho u\otimes u:\nabla\mathcal B\dd x
+\int(\rho^\gamma-\bar\rho^\gamma)(\rho-\bar\rho)\dd x\notag\\
&-2\int\rho^\alpha\nabla u:\nabla\mathcal B\dd x
-2(\alpha-1)\int\rho^\alpha(\divg u)(\rho-\bar\rho)\dd x
-\int\rho^2U^2\dd x.                    \label{eq:cross-energy-identity}
\end{align}
We estimate these terms separately.
First, \eqref{eq:polar-gradient}, the bound $\rho\le M$ and
$\alpha\le1$ give
\begin{equation}\label{eq:kinetic-D0}
\begin{split}
\int\rho|u|^2\dd x
&=\int\rho^{1-\alpha}\rho^\alpha\frac{|u|^2}{r^2}r^2\dd x\\
&\le M^{1-\alpha}R^2\int\rho^\alpha\frac{|u|^2}{r^2}\dd x
\le C D_0(t).
\end{split}
\end{equation}
Consequently,
\[
\left|\int\rho u\otimes u:\nabla\mathcal B\dd x\right|
+\int\rho^2U^2\dd x
\le\bigl(\norm{\nabla\mathcal B}_\infty+M\bigr)
\int\rho|u|^2\dd x\le CD_0(t).
\]
For the pressure term, the pointwise inequality
\[
(s^\gamma-\bar\rho^\gamma)(s-\bar\rho)
\ge\bar\rho^{\gamma-1}(s-\bar\rho)^2,\qquad s\ge0,
\]
implies
\[
\int(\rho^\gamma-\bar\rho^\gamma)(\rho-\bar\rho)\dd x
\ge\bar\rho^{\gamma-1}\norm{\rho-\bar\rho}_2^2.
\]
For the two viscosity terms, Cauchy--Schwarz's inequality,
\eqref{eq:testfield-gradient}, and Young's inequality imply
\begin{align*}
&\left|2\int\rho^\alpha\nabla u:\nabla\mathcal B\dd x
+2(\alpha-1)\int\rho^\alpha(\divg u)(\rho-\bar\rho)\dd x\right|\\
&\qquad\le C D_0(t)^{1/2}\norm{\rho-\bar\rho}_2
\le\frac{\bar\rho^{\gamma-1}}2\norm{\rho-\bar\rho}_2^2
+CD_0(t).
\end{align*}
Substituting these estimates into \eqref{eq:cross-energy-identity},
we find
\begin{equation}\label{eq:cross-energy-lower}
I'(t)\ge\frac{\bar\rho^{\gamma-1}}2\norm{\rho-\bar\rho}_2^2-CD_0(t).
\end{equation}

On the other hand, Poincar\'e's inequality,
\eqref{eq:testfield-gradient} and \eqref{eq:relative-equivalence} yield
\begin{align*}
|I(t)|
&\le\left(\int\rho|u|^2\dd x\right)^{1/2}
\left(\int\rho|\mathcal B|^2\dd x\right)^{1/2}\\
&\le C\left(\int\rho|u|^2\dd x\right)^{1/2}
\norm{\rho-\bar\rho}_2\le C\mathcal E(t).
\end{align*}
Choose a sufficiently small constant $\eta>0$, independent of $T$,
and define $\mathcal E_\eta=\mathcal E-\eta I$. Then
\begin{equation}\label{eq:modified-energy-equivalence}
\frac12\mathcal E(t)\le\mathcal E_\eta(t)\le\frac32\mathcal E(t).
\end{equation}
Combining \eqref{eq:relative-energy-inequality} with
\eqref{eq:cross-energy-lower}, and reducing $\eta$ if necessary, gives
\[
\mathcal E_\eta'(t)+(2\alpha-1)D_0(t)
+\frac{\eta\bar\rho^{\gamma-1}}2\norm{\rho-\bar\rho}_2^2\le0.
\]
By \eqref{eq:relative-equivalence}, \eqref{eq:kinetic-D0} and
\eqref{eq:modified-energy-equivalence},
\[
\mathcal E_\eta(t)\le C\left(D_0(t)+\norm{\rho-\bar\rho}_2^2\right).
\]
Hence there exists a constant $\lambda>0$ such that
\[
\mathcal E_\eta'(t)+\lambda\mathcal E_\eta(t)\le0.
\]
Integrating this differential inequality over $(0,t)$ and using
\eqref{eq:relative-equivalence} and
\eqref{eq:modified-energy-equivalence}, we obtain
\[
\int\rho|u(t)|^2\dd x+\norm{\rho(t)-\bar\rho}_2^2
\le C\mathcal E_\eta(0)e^{-\lambda t}\le Ce^{-\lambda t}.
\]
Finally, integration of \eqref{eq:relative-energy-inequality}
with respect to time over $(t,T)$ yields
\[
2(2\alpha-1)\int_t^TD_0(s)\dd s
\le\mathcal E(t)-\mathcal E(T)\le Ce^{-\lambda t}.
\]
This proves \eqref{eq:relative-decay}.
\end{proof}


\section{The gradient stress for \texorpdfstring{$1/2<\alpha<1$}{1/2 < alpha < 1}}\label{oa:section}

In this section, we take $\sigma=1$ and $1/2<\alpha<1$, and set
\begin{equation}\label{oa:N}
 N_\alpha=\frac{2\alpha\bigl(\alpha+\sqrt{2\alpha-1}\bigr)}
 {(1-\alpha)^2}.
\end{equation}
All the estimates below are established on an arbitrary interval $[0,T]$
on which the local solution exists and the density is positive.
Unless otherwise stated, the constants are independent of $T$ and of
the lower bound of the density on this interval.
We use $D_p$ as defined in \eqref{eq:Dp}, together with the polar
identity \eqref{eq:polar-gradient}.

\subsection{Higher integrability of the velocity}

We start with an elementary algebraic lemma. This lemma gives the precise
range of exponents for which the viscosity controls the gradient after
testing the second equation in \eqref{eq:model} by a power of the velocity.

\Needspace{16\baselineskip}
\begin{lemma}\label{oa:matrix}
Let $1/2<\alpha<1$ and $0\le p<N_\alpha-2$.
There exists a constant $\delta_{\alpha,p}>0$, depending only on
$\alpha$ and $p$, such that, for every real $2\times2$ matrix $A$,
every unit vector $z\in\mathbb R^2$, and every $\tau\in[0,p]$,
\begin{align}
&2A:\bigl((I+\tau z\otimes z)A\bigr)
+2(\alpha-1)\operatorname{tr}A
  \bigl(\operatorname{tr}A+\tau z^{\mathsf T}Az\bigr)
\ge\delta_{\alpha,p}|A|^2.                         \label{oa:matrixineq}
\end{align}
Consequently, for every differentiable vector field $u$,
\begin{align}
&2\nabla u:\nabla(|u|^pu)
+2(\alpha-1)\divg u\,\divg(|u|^pu)
\ge\delta_{\alpha,p}|u|^p|\nabla u|^2.             \label{oa:homogeneous}
\end{align}
When $\alpha=1$, the same statements hold for every finite $p\ge0$.
\end{lemma}

\begin{proof}
We use the convention $A_{ij}=\partial_j u_i$.
Let $O$ be an orthogonal matrix satisfying $Oz=e_1$, and put
$\widetilde A=OAO^{\mathsf T}$.
For any two matrices $A,B$, we have
\[
 (OAO^{\mathsf T}):(OBO^{\mathsf T})
 =\operatorname{tr}\bigl(OA^{\mathsf T}BO^{\mathsf T}\bigr)
 =\operatorname{tr}(A^{\mathsf T}B)=A:B.
\]
Moreover,
\[
 \operatorname{tr}\widetilde A=\operatorname{tr}A,
 \qquad
 e_1^{\mathsf T}\widetilde A e_1=z^{\mathsf T}Az,
 \qquad
 O(I+\tau z\otimes z)AO^{\mathsf T}
 =(I+\tau e_1\otimes e_1)\widetilde A .
\]
Writing
$\widetilde A=\left(\begin{smallmatrix}a&b\\c&d\end{smallmatrix}\right)$,
we therefore have
\begin{align}
&2A:\bigl((I+\tau z\otimes z)A\bigr)
 +2(\alpha-1)\operatorname{tr}A
       \bigl(\operatorname{tr}A+\tau z^{\mathsf T}Az\bigr)\notag\\
&\quad=2\widetilde A:
       \bigl((I+\tau e_1\otimes e_1)\widetilde A\bigr)
 +2(\alpha-1)\operatorname{tr}\widetilde A
       \bigl(\operatorname{tr}\widetilde A
                      +\tau e_1^{\mathsf T}\widetilde A e_1\bigr)\notag\\
&\quad=2\bigl[\alpha(\tau+1)a^2+\alpha d^2
 -(1-\alpha)(\tau+2)ad+(\tau+1)b^2+c^2\bigr],
                                                        \label{oa:quadratic}
\end{align}
and $|A|^2=|\widetilde A|^2=a^2+b^2+c^2+d^2$.
Completing the square gives
\begin{align*}
&\alpha(\tau+1)a^2+\alpha d^2-(1-\alpha)(\tau+2)ad\\
&\quad=\frac{(1-\alpha)(\tau+2)}{2\sqrt{\tau+1}}
       \bigl(\sqrt{\tau+1}\,a-d\bigr)^2\\
&\qquad\quad+\left(\alpha-
       \frac{(1-\alpha)(\tau+2)}{2\sqrt{\tau+1}}\right)
       \bigl((\tau+1)a^2+d^2\bigr).
\end{align*}
Since $(\tau+2)/\sqrt{\tau+1}$ is nondecreasing for $\tau\ge0$ and
$p<N_\alpha-2$, we may take
\[
\delta_{\alpha,p}=2\alpha-
\frac{(1-\alpha)(p+2)}{\sqrt{p+1}}>0.
\]
Dropping the square term and using $0<\delta_{\alpha,p}\le2$ in
\eqref{oa:quadratic} proves \eqref{oa:matrixineq}.

At a point where $u\ne0$, take $z=u/|u|$ and $\tau=p$.
The identities
\[
 \nabla(|u|^pu)=|u|^p(I+pz\otimes z)\nabla u,
 \qquad
 \divg(|u|^pu)=|u|^p\bigl(\divg u+pz^{\mathsf T}(\nabla u)z\bigr)
\]
give \eqref{oa:homogeneous}. At $u=0$, the same inequality follows
directly from the derivative of $u\mapsto |u|^pu$; the case $p=0$
is interpreted as the identity map.
Finally, when $\alpha=1$, the expression in \eqref{oa:quadratic} is
$2[(\tau+1)(a^2+b^2)+c^2+d^2]\ge2|A|^2$, so the statements hold
for every finite $p\ge0$, with $\delta_{1,p}=2$.
This completes the proof.
\end{proof}

We now apply Lemma \ref{oa:matrix} to obtain higher integrability of
the velocity. We treat small positive exponents separately in order to
keep the pressure remainder integrable in time.

\begin{proposition}\label{oa:velocity}
For every $0<p<N_\alpha-2$, we have
\begin{equation}\label{oa:velocityestimate}
 \sup_{0\le t\le T}\int \rho|u|^{p+2}\,dx
 +\int_0^T D_p(t)\,dt\le C_p.
\end{equation}
The constant $C_p$ depends only on the fixed parameters, $p$, and the
initial data.
\end{proposition}

\begin{proof}
We first consider $0<p\le2$. We follow Lei and Liang \cite{LL2026}
in combining a regularized velocity test with the internal-energy
identity. Put $\Phi=(1+|u|^2)^{p/2}$. Multiplying the second equation in
\eqref{eq:model} by $\Phi u$ and using the first equation, we obtain
\[
\frac1{p+2}\frac{d}{dt}\int\rho
 \bigl((1+|u|^2)^{(p+2)/2}-1\bigr)\,dx
+\mathcal D_\Phi
=\int\rho^\gamma\divg(\Phi u)\,dx,
\]
where
\begin{align*}
\mathcal D_\Phi={}&
 2\int\rho^\alpha\nabla u:\nabla(\Phi u)\,dx
 +2(\alpha-1)\int\rho^\alpha\divg u\,\divg(\Phi u)\,dx.
\end{align*}
On the other hand, the first equation in \eqref{eq:model} gives
\[
\frac{d}{dt}\int\frac{\rho^\gamma}{\gamma-1}\,dx
=-\int\rho^\gamma\divg u\,dx.
\]
Adding these two identities, we obtain
\begin{equation}\label{oa:regularentropy}
\begin{split}
&\frac{d}{dt}\int
 \left\{\frac{\rho}{p+2}
 \bigl((1+|u|^2)^{(p+2)/2}-1\bigr)
 +\frac{\rho^\gamma}{\gamma-1}\right\}\,dx
 +\mathcal D_\Phi\\
&\qquad=\int\rho^\gamma\divg\bigl((\Phi-1)u\bigr)\,dx.
\end{split}
\end{equation}

At $u\ne0$, we have
\[
 \nabla(\Phi u)=\Phi(I+\tau z\otimes z)\nabla u,
 \qquad z=\frac{u}{|u|},\qquad
 \tau=\frac{p|u|^2}{1+|u|^2}\in[0,p].
\]
Lemma \ref{oa:matrix}, including its direct interpretation at $u=0$,
therefore gives
\begin{equation}\label{oa:regularcoercivity}
 \mathcal D_\Phi\ge
 \delta_{\alpha,p}\int 
 \rho^\alpha(1+|u|^2)^{p/2}|\nabla u|^2\,dx .
\end{equation}
For $x\ge0$ and $0<p\le2$,
\[
 (1+x)^{p/2}-1\le x(1+x)^{p/2-1}.
\]
Differentiating $(\Phi-1)u$ consequently yields
\[
 \bigl|\divg((\Phi-1)u)\bigr|
 \le C_p|u|^2(1+|u|^2)^{p/2-1}|\nabla u|.
\]
By Young's inequality, the right-hand side of
\eqref{oa:regularentropy} is at most
\begin{align*}
 &\frac{\delta_{\alpha,p}}2
  \int \rho^\alpha(1+|u|^2)^{p/2}|\nabla u|^2\,dx\\
 &\quad
 +C\int \rho^{2\gamma-\alpha}
 |u|^4(1+|u|^2)^{p/2-2}\,dx\\
 &\le \frac{\delta_{\alpha,p}}2
  \int \rho^\alpha(1+|u|^2)^{p/2}|\nabla u|^2\,dx
 +C\int \rho^{2\gamma-\alpha}|u|^2\,dx\\
 &\le \frac{\delta_{\alpha,p}}2
  \int \rho^\alpha(1+|u|^2)^{p/2}|\nabla u|^2\,dx
 +C\norm{\rho^{2\gamma-2\alpha}r^2}_{L^\infty}D_0(t).
\end{align*}
Here $|u|^4(1+|u|^2)^{p/2-2}\le|u|^2$, and we used
$|\nabla u|^2\ge |u|^2/r^2$ in the last step.
Since $2\gamma-2\alpha>0$, \eqref{eq:sphere} bounds the coefficient.
Combining this estimate with \eqref{oa:regularentropy} and
\eqref{oa:regularcoercivity}, we obtain
\begin{equation}\label{oa:smallp-combined}
\begin{split}
&\frac{d}{dt}\int
 \left\{\frac{\rho}{p+2}
 \bigl((1+|u|^2)^{(p+2)/2}-1\bigr)
 +\frac{\rho^\gamma}{\gamma-1}\right\}\,dx\\
&\qquad+\frac{\delta_{\alpha,p}}2
 \int\rho^\alpha(1+|u|^2)^{p/2}|\nabla u|^2\,dx
 \le CD_0(t).
\end{split}
\end{equation}
Integrating with respect to time and using
$(1+|u|^2)^{(p+2)/2}-1\ge|u|^{p+2}$ and \eqref{eq:basic}, we find
\begin{align*}
&\sup_{0\le t\le T}\int\rho|u|^{p+2}\,dx
 +\int_0^TD_p(t)\,dt\\
&\qquad\le C\int
 \left(\rho_0|u_0|^{p+2}+\rho_0|u_0|^2+\rho_0^\gamma\right)\,dx
 +C\int_0^TD_0(t)\,dt\le C_p.
\end{align*}

If $2<p<N_\alpha-2$, we multiply the second equation in
\eqref{eq:model} by $|u|^pu$.
Lemma \ref{oa:matrix} and Young's inequality give
\begin{align}
 \frac1{p+2}\frac{d}{dt}\int \rho|u|^{p+2}\,dx
 +\frac{\delta_{\alpha,p}}2D_p(t)
 \le C\int \rho^{2\gamma-\alpha}|u|^p\,dx.
                                                        \label{oa:largep}
\end{align}
The last term can be estimated as follows:
\begin{align*}
 \int \rho^{2\gamma-\alpha}|u|^p\,dx
 &=\int 
 \left(\rho^\alpha\frac{|u|^{p+2}}{r^2}\right)^{(p-2)/p}
 \left(\rho^\alpha\frac{|u|^2}{r^2}\right)^{2/p}
 \rho^{2\gamma-2\alpha}r^2\,dx\\
 &\le\varepsilon D_p(t)
 +C_{\varepsilon,p}
 \int \rho^\alpha\frac{|u|^2}{r^2}
       \bigl(\rho^{2\gamma-2\alpha}r^2\bigr)^{p/2}\,dx\\
 &\le\varepsilon D_p(t)+C_{\varepsilon,p}D_0(t).
\end{align*}
We again used \eqref{eq:sphere}.
Substituting this estimate into \eqref{oa:largep} and choosing
$\varepsilon$ sufficiently small, we obtain
\[
\frac1{p+2}\frac{d}{dt}\int\rho|u|^{p+2}\,dx
+\frac{\delta_{\alpha,p}}4D_p(t)\le CD_0(t).
\]
Integrating with respect to time and using \eqref{eq:basic} yields
\[
\sup_{0\le t\le T}\int\rho|u|^{p+2}\,dx
+\int_0^TD_p(t)\,dt
\le C\int\rho_0|u_0|^{p+2}\,dx
+C\int_0^TD_0(t)\,dt\le C_p.
\]
This completes the proof.
\end{proof}

\subsection{The negative part of the effective velocity}

The radial component of the effective velocity is
\begin{equation}\label{oa:Vdefinition}
 V=U+2\alpha\rho^{\alpha-2}\rho_r .
\end{equation}
It satisfies
\begin{equation}\label{oa:effective}
 \rho(V_t+UV_r)
 +\frac{\gamma}{2\alpha}\rho^{\gamma-\alpha+1}V
 =\frac{\gamma}{2\alpha}\rho^{\gamma-\alpha+1}U
 +\rho\frac{W^2}{r}.
\end{equation}
The sign of the last term allows us to estimate the negative part of
$V$ before controlling the full effective velocity.

\begin{lemma}\label{oa:negativeV}
For every $2<m<N_\alpha$, there holds
\begin{equation}\label{oa:negativeVestimate}
 \sup_{0\le t\le T}\int \rho V_-^m\,dx\le C_m,
 \qquad V_-=\max\{-V,0\}.
\end{equation}
Moreover, there is a constant $K\ge1$ such that
\begin{equation}\label{oa:onesided}
 \rho(r,t)\le K\rho(s,t),
 \qquad 0\le r\le s\le R,\quad 0\le t\le T,
\end{equation}
and
\begin{equation}\label{oa:upperrho}
 \sup_{0\le t\le T}\norm{\rho(t)}_{L^\infty}\le C .
\end{equation}
\end{lemma}

\begin{proof}
Multiplying \eqref{oa:effective} by $-V_-^{m-1}$, integrating over
$\Omega$, and using the first equation in \eqref{eq:model}, we obtain
\begin{align*}
 &\frac1m\frac{d}{dt}\int \rho V_-^m\,dx
 +\frac{\gamma}{2\alpha}
       \int \rho^{\gamma-\alpha+1}V_-^m\,dx\\
 &\quad=-\frac{\gamma}{2\alpha}
       \int \rho^{\gamma-\alpha+1}U V_-^{m-1}\,dx
       -\int \frac{\rho W^2}{r}V_-^{m-1}\,dx\\
 &\quad\le\frac{\gamma}{4\alpha}
       \int \rho^{\gamma-\alpha+1}V_-^m\,dx
       +C_m\int \rho^{\gamma-\alpha+1}|U|^m\,dx.
\end{align*}
Here the last term on the second line is nonpositive, and the remaining
term is estimated by Young's inequality.
Since $\gamma-2\alpha+1>0$, \eqref{eq:sphere} gives
\begin{equation}\label{oa:pressureU}
 \int \rho^{\gamma-\alpha+1}|U|^m\,dx
 \le\norm{\rho^{\gamma-2\alpha+1}r^2}_{L^\infty}
 \int \rho^\alpha\frac{|u|^m}{r^2}\,dx
 \le C D_{m-2}(t).
\end{equation}
Consequently, integrating with respect to time over $(0,t)$ gives
\[
\int \rho V_-^m(t)\,dx
+\frac{m\gamma}{4\alpha}\int_0^t\int 
       \rho^{\gamma-\alpha+1}V_-^m\,dx\,ds
\le\int \rho_0 V_-(0)^m\,dx
+C_m\int_0^tD_{m-2}(s)\,ds\le C_m,
\]
where the last inequality follows from Proposition \ref{oa:velocity}.
Taking the supremum for $0\le t\le T$ proves
\eqref{oa:negativeVestimate}.

We next take $m_c=1/(1-\alpha)>2$.
This exponent belongs to $(2,N_\alpha)$. To see this, it is enough to
notice that
\[
 2\alpha\bigl(\alpha+\sqrt{2\alpha-1}\bigr)>1-\alpha,
 \qquad
 2\alpha^2-(1-\alpha)=(2\alpha-1)(\alpha+1)>0.
\]
It follows from \eqref{oa:Vdefinition} that
\[
 \bigl((\log\rho)_r\bigr)_-
 \le\frac1{2\alpha}\rho^{1-\alpha}(V_-+|U|).
\]
As $m_c(1-\alpha)=1$, Proposition \ref{oa:velocity} and
\eqref{oa:negativeVestimate} imply
\[
 \sup_{0\le t\le T}\int_0^R
 \bigl|((\log\rho)_r)_-\bigr|^{m_c}r\,dr\le C .
\]
For $0\le r\le s\le R$, radial H\"older's inequality gives
\begin{align*}
 \log\rho(r,t)-\log\rho(s,t)
 &\le\int_r^s\bigl((\log\rho)_\xi\bigr)_-\,d\xi\\
 &\le C\left(\int_0^R
     \xi^{-1/(m_c-1)}\,d\xi\right)^{(m_c-1)/m_c}
 \le C .
\end{align*}
The integral is finite because $m_c>2$.
Exponentiating this inequality proves \eqref{oa:onesided}.

Finally, Lemma \ref{lem:basic} gives $\rho(R,t)\le C$.
Taking $s=R$ in \eqref{oa:onesided} therefore proves
\eqref{oa:upperrho}. This completes the proof.
\end{proof}

\subsection{Weighted estimates for the tangential velocity}

\begin{lemma}\label{oa:allangular}
For every $q\ge2$ and $0\le k<2\sqrt{q-1}$, we have
\begin{align}
 &\sup_{0\le t\le T}\int\rho|W|^qr^{-k}\,dx\notag\\
 &\quad+\int_0^T\int\rho^\alpha|W|^{q-2}r^{-k}
       \left(W_r^2+\frac{W^2}{r^2}\right)\,dx\,dt
 \le C_{q,k}.                                          \label{oa:allangularestimate}
\end{align}
\end{lemma}

\begin{proof}
Multiplying \eqref{eq:angular} by $|W|^{q-2}Wr^{-k}$ and using
the first equation in \eqref{eq:model}, we obtain
\begin{align}
 &\frac1q\frac{d}{dt}\int\rho|W|^qr^{-k}\,dx
 +2\int\rho^\alpha|W|^{q-2}r^{-k}
 \left((q-1)W_r^2-kW_r\frac Wr+\frac{W^2}{r^2}\right)\,dx\notag\\
 &\qquad=-\left(1+\frac{k}{q}\right)
             \int\rho U|W|^qr^{-k-1}\,dx.          \label{oa:angularidentity}
\end{align}
Here integration by parts gives
\begin{align*}
 &4\pi\int_0^R(r\rho^\alpha W_r)_r
               |W|^{q-2}Wr^{-k}\,dr
 -4\pi\int_0^R\rho^\alpha|W|^qr^{-k-1}\,dr\\
 &\qquad=-4\pi\int_0^R\rho^\alpha|W|^{q-2}r^{1-k}
 \left((q-1)W_r^2-kW_r\frac Wr+\frac{W^2}{r^2}\right)\,dr.
\end{align*}

Since $k<2\sqrt{q-1}$,
\[
(q-1)W_r^2-kW_r\frac Wr+\frac{W^2}{r^2}
\ge\left(1-\frac{k}{2\sqrt{q-1}}\right)
       \left(W_r^2+\frac{W^2}{r^2}\right)
=:c_{q,k}\left(W_r^2+\frac{W^2}{r^2}\right),
\]
where $c_{q,k}>0$.

We next consider the right-hand side of \eqref{oa:angularidentity}. To this end, we choose 
$\frac{2\alpha}{1-\alpha}<S<N_\alpha.$ This is possible since $\alpha+\sqrt{2\alpha-1}>1-\alpha$.
By \eqref{oa:onesided} and $u(R,t)=0$,
\begin{align}
 \rho(r,t)^\alpha|u(r,t)|^S
 &\le S\rho(r,t)^\alpha
       \int_r^R|u(s,t)|^{S-1}|\partial_su(s,t)|\,ds\notag\\
 &\le C\int_r^R\rho(s,t)^\alpha
       \left(\frac{|u|^S}{s}
                   +s|u|^{S-2}|\partial_su|^2\right)\,ds
 \le CD_{S-2}(t).                                  \label{oa:pointwiseu}
\end{align}
Consequently, \eqref{oa:upperrho} and $1-\alpha-2\alpha/S>0$ give
\begin{equation}\label{oa:coefficient}
 \rho^{2-\alpha}U^2
 \le CD_{S-2}(t)^{2/S}\rho^{2-\alpha-2\alpha/S}
 \le CD_{S-2}(t)^{2/S}\rho.
\end{equation}
Moreover, \eqref{oa:upperrho} implies
\begin{equation}\label{oa:BcontrolsZ}
\int\rho|W|^qr^{-k}\,dx
\le R^2\norm{\rho}_\infty^{1-\alpha}
       \int\rho^\alpha|W|^qr^{-k-2}\,dx
\le C\int\rho^\alpha|W|^qr^{-k-2}\,dx.
\end{equation}
By Young's inequality and \eqref{oa:coefficient},
\begin{align*}
&\left(1+\frac{k}{q}\right)
 \int\rho|U||W|^qr^{-k-1}\,dx\\
&\quad\le\varepsilon\int\rho^\alpha|W|^qr^{-k-2}\,dx
       +C_\varepsilon D_{S-2}(t)^{2/S}
                      \int\rho|W|^qr^{-k}\,dx\\
&\quad\le\varepsilon\int\rho^\alpha|W|^qr^{-k-2}\,dx
       +\eta\int\rho|W|^qr^{-k}\,dx
       +C_{\varepsilon,\eta}D_{S-2}(t)
                      \int\rho|W|^qr^{-k}\,dx\\
&\quad\le(\varepsilon+C\eta)
                      \int\rho^\alpha|W|^qr^{-k-2}\,dx
       +C_{\varepsilon,\eta}D_{S-2}(t)
                      \int\rho|W|^qr^{-k}\,dx.
\end{align*}
Here $2/S<1$ and we used \eqref{oa:BcontrolsZ} in the last step.
Choosing $\varepsilon,\eta>0$ such that
$\varepsilon+C\eta\le c_{q,k}$, we deduce from
\eqref{oa:angularidentity} that
\begin{align}
 &\frac1q\frac{d}{dt}\int\rho|W|^qr^{-k}\,dx
 +c_{q,k}\int\rho^\alpha|W|^{q-2}r^{-k}
             \left(W_r^2+\frac{W^2}{r^2}\right)\,dx\notag\\
 &\qquad\le CD_{S-2}(t)\int\rho|W|^qr^{-k}\,dx.
                                                   \label{oa:angular-differential}
\end{align}
By $k<2\sqrt{q-1}\le q$ and \eqref{eq:polar-gradient},
\[
\int\rho_0|W_0|^qr^{-k}\,dx
\le R^{q-k}\int\rho_0\left|\frac{W_0}{r}\right|^q\,dx
\le R^{q-k}\int\rho_0|\nabla u_0|^q\,dx\le C.
\]
Applying Gronwall's inequality to \eqref{oa:angular-differential} and
using Proposition \ref{oa:velocity}, we conclude that
\begin{align*}
 &\sup_{0\le t\le T}\int\rho|W|^qr^{-k}\,dx
 +\int_0^T\int\rho^\alpha|W|^{q-2}r^{-k}
             \left(W_r^2+\frac{W^2}{r^2}\right)\,dx\,dt\\
 &\qquad\le C\exp\left(C\int_0^TD_{S-2}(t)\,dt\right)
                      \int\rho_0|W_0|^qr^{-k}\,dx
 \le C_{q,k}.
\end{align*}
This completes the proof.
\end{proof}

In particular, taking $k=0$, we have, for every finite $q\ge2$,
\begin{equation}\label{oa:unweightedangular}
 \int_0^T D_q^{\rm ang}(t)\,dt\le C_q,\qquad
 D_q^{\rm ang}(t)=\int\rho^\alpha|W|^{q-2}
       \left(W_r^2+\frac{W^2}{r^2}\right)\,dx.
\end{equation}

\subsection{Higher integrability of the radial effective velocity}

The weighted angular estimates give the following bound for $V$.
We use the notation
\begin{equation}\label{oa:F}
 F=\frac{W^2}{r},\qquad J_M(t)=\int\rho F^M\,dx.
\end{equation}

\begin{proposition}\label{oa:first-effective}
For every $2<m<\min\{N_\alpha,6+2\sqrt7\}$, there exists a constant
$C_m>0$, depending on the fixed parameters, $m$ and the initial data,
but independent of $T$, such that
\begin{equation}\label{oa:firstVmoment}
 \sup_{0\le t\le T}\left(\int\rho|V(t)|^m\,dx\right)^{1/m}
 +\int_0^T\left(\int\rho\frac{|W(t)|^{2m}}{r^m}\,dx\right)^{1/m}\,dt
 \le C_m.
\end{equation}
\end{proposition}

\begin{proof}
For every $2<m<N_\alpha$, multiplying \eqref{oa:effective} by
$|V|^{m-2}V$ and using the first
equation in \eqref{eq:model}, Young's inequality and
\eqref{oa:pressureU}, we obtain
\begin{align}
 &\frac1m\frac{d}{dt}\int\rho|V|^m\,dx
  +c\int\rho^{\gamma-\alpha+1}|V|^m\,dx\notag\\
 &\qquad\le CD_{m-2}(t)
  +\left(\int\rho F^m\,dx\right)^{1/m}
     \left(\int\rho|V|^m\,dx\right)^{(m-1)/m}.
                                                        \label{oa:Vmoment}
\end{align}
We next estimate the centrifugal term.
By Lemma \ref{lem:decay} and \eqref{oa:upperrho},
\begin{equation}\label{oa:energydecay}
 \int\rho|u|^2\,dx+\norm{\rho-\bar\rho}_{L^2}^2
 +\int_t^TD_0(s)\,ds\le Ce^{-\lambda t}.
\end{equation}
Taking $q=2$ and $k=3/2$ in Lemma \ref{oa:allangular}, we obtain
\begin{align}
 J_1(t)
 &\le\left(\int\rho W^2\,dx\right)^{1/3}
       \left(\int\rho W^2r^{-3/2}\,dx\right)^{2/3}
 \le Ce^{-\lambda t/3}.                              \label{oa:J1}
\end{align}
Choose $m<M<6+2\sqrt7$. Taking $q=2M$ and $k=M-2$ in
Lemma \ref{oa:allangular} gives
\begin{equation}\label{oa:firstJM}
 \int_0^T J_M(t)\,dt\le C_M.
\end{equation}
Here we used
\[
 (M-2)^2<4(2M-1),\qquad
 \rho F^M\le C\rho^\alpha|W|^{2M}r^{-M}.
\]
H\"older's inequality yields
\begin{equation}\label{oa:sourceinterpolation}
 \left(\int\rho F^m\,dx\right)^{1/m}
 \le J_1^{(M-m)/(m(M-1))}
      J_M^{(m-1)/(m(M-1))}.
\end{equation}
By \eqref{oa:J1}, \eqref{oa:firstJM} and H\"older's inequality in time,
\begin{align*}
 \int_0^T\left(\int\rho F^m\,dx\right)^{1/m}\,dt
 &\le C\int_0^Te^{-ct}
             J_M(t)^{(m-1)/(m(M-1))}\,dt\\
 &\le C\left(\int_0^TJ_M(t)\,dt\right)^{(m-1)/(m(M-1))}
 \le C_m.
\end{align*}

For every $2<m<N_\alpha$, multiplying \eqref{oa:Vmoment} by
$\left(1+\int\rho|V|^m\,dx\right)^{-(m-1)/m}$ gives
\[
\frac{d}{dt}\left(1+\int\rho|V|^m\,dx\right)^{1/m}
\le CD_{m-2}(t)+C\left(\int\rho F^m\,dx\right)^{1/m}.
\]
Integrating with respect to time, we obtain
\begin{equation}\label{oa:Vroot}
\begin{split}
 &\sup_{0\le t\le T}
       \left(1+\int\rho(t)|V(t)|^m\,dx\right)^{1/m}\\
 &\quad\le\left(1+\int\rho_0|V(0)|^m\,dx\right)^{1/m}
       +C\int_0^TD_{m-2}(t)\,dt\\
 &\qquad+C\int_0^T\left(\int\rho F^m\,dx\right)^{1/m}\,dt.
\end{split}
\end{equation}
The right-hand side is bounded by Proposition \ref{oa:velocity} and
the estimate above. This proves \eqref{oa:firstVmoment}.
\end{proof}

\subsection{Integrated coercivity near the origin}

The pointwise condition $k^2<4(q-1)$ is too restrictive for large
singular weights. Near the origin, it can be replaced by an integrated
inequality under a relative continuity condition on the viscosity.

\begin{lemma}\label{oa:frozen}
Let $q>2$ and $0<k<q$.
There are constants $\varepsilon_{q,k}>0$ and $c_{q,k}>0$ with the
following property. If $0<\delta\le R$ and
\begin{equation}\label{oa:freezingcondition}
 \left|\frac{\rho(r,t)^\alpha}{\rho(0,t)^\alpha}-1\right|
 \le\varepsilon_{q,k},\qquad 0\le r\le\delta,
\end{equation}
then
\begin{align}
 &\int_{|x|<\delta}\rho^\alpha|W|^{q-2}r^{-k}
 \left((q-1)W_r^2-kW_r\frac Wr+\frac{W^2}{r^2}\right)\,dx\notag\\
 &\qquad\ge c_{q,k}\int_{|x|<\delta}\rho^\alpha|W|^{q-2}r^{-k}
 \left(W_r^2+\frac{W^2}{r^2}\right)\,dx .             \label{oa:frozencoercivity}
\end{align}
The constants are independent of $\rho(0,t)$ and $\delta$.
No boundary condition at $r=\delta$ is required.
\end{lemma}

\begin{proof}
We first take the coefficient $\rho^\alpha$ to be constant.
Omitting the common angular factor, put
\[
 G=\int_0^\delta |W|^{q-2}W_r^2r^{1-k}\,dr,\quad
 B=\int_0^\delta |W|^qr^{-k-1}\,dr,\quad
 L=\delta^{-k}|W(\delta)|^q\ge0 .
\]
Integration by parts gives
\begin{equation}\label{oa:crossboundary}
 \int_0^\delta |W|^{q-2}WW_rr^{-k}\,dr
 =\frac{kB+L}{q}.
\end{equation}
The origin contributes zero because $W=O(r)$ and $q>k$.
By completing the square,
\begin{align*}
 0&\le\int_0^\delta |W|^{q-2}
       \left(W_r-\frac{kW}{qr}\right)^2r^{1-k}\,dr\\
 &=G-\frac{2k}{q}\frac{kB+L}{q}+\frac{k^2}{q^2}B ,
\end{align*}
and hence
\begin{equation}\label{oa:Hardyboundary}
 G\ge\frac{k^2}{q^2}B+\frac{2k}{q^2}L .
\end{equation}
Choose $\varepsilon>0$ sufficiently small that
$q-2-2\varepsilon>0$ and $(1+\varepsilon)k^2<q^2$.
Using \eqref{oa:crossboundary} and
\eqref{oa:Hardyboundary}, we obtain
\begin{align*}
 (q-1)G-\frac{k(kB+L)}q+B
 &\ge \varepsilon G
  +\left(1-\frac{(1+\varepsilon)k^2}{q^2}\right)B
  +\frac{k(q-2-2\varepsilon)}{q^2}L\\
 &\ge c^{(0)}_{q,k}(G+B),
\end{align*}
where $c^{(0)}_{q,k}>0$ depends only on $q,k$.

For the variable coefficient, take
$\varepsilon_{q,k}=\min\{1/2,c^{(0)}_{q,k}/(2q-2+k)\}$.
By Young's inequality and $q>2$,
\begin{align*}
\left|(q-1)W_r^2-kW_r\frac Wr+\frac{W^2}{r^2}\right|
&\le(q-1)W_r^2+k|W_r|\frac{|W|}{r}+\frac{W^2}{r^2}\\
&\le\left(q-1+\frac k2\right)
       \left(W_r^2+\frac{W^2}{r^2}\right).
\end{align*}
Splitting off the constant coefficient and using
\eqref{oa:freezingcondition}, we obtain
\begin{align*}
&\int_0^\delta\rho^\alpha|W|^{q-2}r^{1-k}
  \left((q-1)W_r^2-kW_r\frac Wr+\frac{W^2}{r^2}\right)\,dr\\
&\quad=\rho(0,t)^\alpha
       \left((q-1)G-\frac{k(kB+L)}q+B\right)\\
&\qquad\quad+\int_0^\delta
  \bigl(\rho^\alpha-\rho(0,t)^\alpha\bigr)|W|^{q-2}r^{1-k}
  \left((q-1)W_r^2-kW_r\frac Wr+\frac{W^2}{r^2}\right)\,dr\\
&\quad\ge\rho(0,t)^\alpha
 \left[(q-1)G-\frac{k(kB+L)}q+B
       -\varepsilon_{q,k}\left(q-1+\frac k2\right)(G+B)\right]\\
&\quad\ge\frac{c^{(0)}_{q,k}}2\rho(0,t)^\alpha(G+B)\\
&\quad\ge\frac{c^{(0)}_{q,k}}3
 \int_0^\delta\rho^\alpha|W|^{q-2}r^{1-k}
             \left(W_r^2+\frac{W^2}{r^2}\right)\,dr.
\end{align*}
The last step uses $\rho^\alpha\le(1+\varepsilon_{q,k})\rho(0,t)^\alpha
\le3\rho(0,t)^\alpha/2$.
Multiplying by $2\pi$ proves \eqref{oa:frozencoercivity}, with
$c_{q,k}=c^{(0)}_{q,k}/3$.
\end{proof}

\subsection{A conditional estimate with a precise loss}

On the fixed local existence interval, let
\begin{equation}\label{oa:H}
 h(t)=\max\{1,\rho(0,t)^{-1}\},\qquad
 H=\max_{0\le t\le T}h(t).
\end{equation}
Local positivity ensures that $H<\infty$.
We will first keep track of its power in the estimates, and then
bound $H$ by a finite iteration.

\begin{lemma}\label{oa:conditional}
Fix $M>m>2$ with $m<N_\alpha$, and let
\begin{equation}\label{oa:conditionalindices}
 A=\frac{(M-2)(m-1)}{m(M-1)}.
\end{equation}
Suppose that, at each time $t\in[0,T]$, there is a radius
$\delta(t)\in(0,R/2]$ such that
\[
\left|\frac{\rho(r,t)^\alpha}{\rho(0,t)^\alpha}-1\right|
\le\varepsilon_{2M,M-2},\qquad 0\le r\le\delta(t),
\]
and
\begin{equation}\label{oa:conditionalradius}
 \delta(t)\ge c_dH^{-d},
\end{equation}
where $d\ge0$ and $c_d>0$ are independent of $H$ and $T$.
Then
\begin{align}
 &\int_0^TJ_M(t)\,dt\le C_d H^{d(M-2)},           \label{oa:conditionalJM}\\
 &\sup_{0\le t\le T}\left(\int\rho|V(t)|^m\,dx\right)^{1/m}
 \le C_d H^{Ad}.                               \label{oa:conditionalV}
\end{align}
\end{lemma}

\begin{proof}
Choose $2\alpha/(1-\alpha)<S<N_\alpha$.
Testing \eqref{eq:angular} by $|W|^{2M-2}Wr^{2-M}$ and applying
Lemma \ref{oa:frozen} on $|x|<\delta(t)$, we obtain
\begin{align*}
&\frac{d}{dt}\int\rho|W|^{2M}r^{2-M}\,dx
 +c\int\rho^\alpha|W|^{2M}r^{-M}\,dx\\
&\qquad\le C\int\rho|U||W|^{2M}r^{1-M}\,dx
 +C\delta(t)^{2-M}D_{2M}^{\rm ang}(t).
\end{align*}
Here the outer part of the quadratic integral and the outer positive
term are estimated by
\begin{align*}
&\left|\int_{|x|\ge\delta(t)}\rho^\alpha|W|^{2M-2}r^{2-M}
 \left((2M-1)W_r^2-(M-2)W_r\frac Wr+\frac{W^2}{r^2}\right)\,dx\right|\\
&\quad+\int_{|x|\ge\delta(t)}\rho^\alpha|W|^{2M}r^{-M}\,dx\\
&\qquad\le C\int_{|x|\ge\delta(t)}\rho^\alpha|W|^{2M-2}r^{2-M}
             \left(W_r^2+\frac{W^2}{r^2}\right)\,dx
 \le C\delta(t)^{2-M}D_{2M}^{\rm ang}(t).
\end{align*}
Using \eqref{oa:coefficient} and \eqref{oa:BcontrolsZ} as in
Lemma \ref{oa:allangular}, we have
\begin{align*}
\int\rho|U||W|^{2M}r^{1-M}\,dx
&\le\varepsilon\int\rho^\alpha|W|^{2M}r^{-M}\,dx\\
&\quad+C_\varepsilon D_{S-2}(t)
                     \int\rho|W|^{2M}r^{2-M}\,dx.
\end{align*}
Taking $\varepsilon$ sufficiently small and using
\eqref{oa:conditionalradius}, we deduce that
\begin{equation}\label{oa:conditionalZ}
\begin{split}
&\frac{d}{dt}\int\rho|W|^{2M}r^{2-M}\,dx
 +c\int\rho^\alpha|W|^{2M}r^{-M}\,dx\\
&\qquad\le CD_{S-2}(t)\int\rho|W|^{2M}r^{2-M}\,dx
 +C_dH^{d(M-2)}D_{2M}^{\rm ang}(t).
\end{split}
\end{equation}
Since
\[
\int\rho_0|W_0|^{2M}r^{2-M}\,dx
\le R^{M+2}\int\rho_0|\nabla u_0|^{2M}\,dx\le C,
\]
Gronwall's inequality, Proposition \ref{oa:velocity} and
\eqref{oa:unweightedangular} give
\[
\int_0^T\int\rho^\alpha|W|^{2M}r^{-M}\,dx\,dt
\le C_dH^{d(M-2)}.
\]
Together with $J_M(t)\le C\int\rho^\alpha|W|^{2M}r^{-M}\,dx$,
this proves \eqref{oa:conditionalJM}.

Finally, \eqref{oa:sourceinterpolation}, \eqref{oa:J1} and H\"older's
inequality in time imply
\begin{align*}
\int_0^T\left(\int\rho F^m\,dx\right)^{1/m}\,dt
&\le C\left(\int_0^TJ_M(t)\,dt\right)^{\frac{m-1}{m(M-1)}}\\
&\le C_dH^{\frac{d(M-2)(m-1)}{m(M-1)}}
 =C_dH^{Ad}.
\end{align*}
Substituting this estimate into \eqref{oa:Vroot} and using
Proposition \ref{oa:velocity} proves \eqref{oa:conditionalV}.
\end{proof}

\subsection{A finite iteration for the lower density bound}

We now combine the conditional estimate in Lemma \ref{oa:conditional}
with two estimates of the freezing radius. A finite iteration will
remove the dependence on the inverse central density $H$.

\begin{proposition}\label{oa:lower}
There exists a constant $C>0$, independent of $T$, such that
\begin{equation}\label{oa:lowerestimate}
 C^{-1}\le \rho(r,t)\le C,\qquad
 0\le r\le R,\quad 0\le t\le T .
\end{equation}
\end{proposition}

\begin{proof}
The upper bound has already been proved in Lemma \ref{oa:negativeV}.
To prove the lower bound, choose
\begin{align}
 &2<m_0<\min\left\{\frac1{1-\alpha},N_\alpha,6+2\sqrt7\right\},
                                                        \label{oa:m0}\\
 &\frac1{1-\alpha}<m<N_\alpha.                           \label{oa:mhigh}
\end{align}
Let
\begin{equation}\label{oa:initialindices}
 d_0=\frac{m_0(\alpha-1)+1}{m_0-2}>0,
\end{equation}
and choose $M>m$ sufficiently close to $m$ that
\begin{equation}\label{oa:contractionchoice}
 \frac{M-m}{M-1}\,d_0<m(1-\alpha)-1.
\end{equation}
The constant $A$ is given by \eqref{oa:conditionalindices}.
All these exponents are fixed in what follows.

By \eqref{oa:Vdefinition}, Propositions \ref{oa:velocity} and
\ref{oa:first-effective}, and radial H\"older's inequality,
\begin{align}
 &\left|\rho(r,t)^{\alpha-1+1/m_0}
          -\rho(0,t)^{\alpha-1+1/m_0}\right|\notag\\
 &\qquad\le C\left(\int_0^r\rho|V-U|^{m_0}s\,ds\right)^{1/m_0}
                r^{1-2/m_0}
 \le Cr^{1-2/m_0}.                                  \label{oa:positiveholder}
\end{align}
Since $\rho(0,t)^{\alpha-1+1/m_0}\ge h(t)^{-(\alpha-1+1/m_0)}$,
a sufficiently small constant $c_0>0$ gives a radius
\begin{equation}\label{oa:baseline}
 \delta_0(t)=c_0h(t)^{-d_0}\le R/2,\qquad
 \delta_0(t)\ge c_0H^{-d_0},
\end{equation}
on which
\[
 \left|\frac{\rho(r,t)^\alpha}{\rho(0,t)^\alpha}-1\right|
 \le\varepsilon_{2M,M-2}.
\]
Indeed, \eqref{oa:positiveholder} and the choice of $d_0$ imply
\[
\frac{\left|\rho(r,t)^{\alpha-1+1/m_0}
              -\rho(0,t)^{\alpha-1+1/m_0}\right|}
     {\rho(0,t)^{\alpha-1+1/m_0}}
\le Cc_0^{1-2/m_0},\qquad r\le\delta_0(t).
\]

We next improve the power of $H$ in the radius estimate.
Suppose that a radius satisfying the same relative-density condition
has been obtained with $\delta(t)\ge c_dH^{-d}$.
Lemma \ref{oa:conditional}, Proposition \ref{oa:velocity} and
\eqref{oa:Vdefinition} give
\begin{align}
 &\left|\rho(r,t)^{\alpha-1+1/m}
          -\rho(0,t)^{\alpha-1+1/m}\right|\notag\\
 &\qquad\le C\left(\int_0^r\rho|V-U|^m s\,ds\right)^{1/m}
                    r^{1-2/m}
 \le C_dH^{Ad}r^{1-2/m}.                            \label{oa:negativeholder}
\end{align}
The density upper bound implies
$\rho(0,t)^{\alpha-1+1/m}\ge C^{-1}h(t)^{1-\alpha-1/m}$.
Thus the relative-density condition also holds on
\begin{equation}\label{oa:secondradius}
 \delta_1(t)=\min\left\{\frac R2,
 c'_dH^{-\frac{Adm}{m-2}}
       h(t)^{\frac{m(1-\alpha)-1}{m-2}}\right\},
\end{equation}
provided $c'_d>0$ is sufficiently small.
Both radii satisfy the same condition, so their maximum does as well.
If $\delta_1(t)<R/2$, then
\begin{align}
 &\max\{\delta_0(t),\delta_1(t)\}^{(m-2)d_0+m(1-\alpha)-1}\notag\\
 &\qquad\ge\delta_0(t)^{m(1-\alpha)-1}
             \delta_1(t)^{(m-2)d_0}
 \ge C_dH^{-Adm d_0}.                               \label{oa:combine}
\end{align}
Consequently, a valid radius is bounded below by
$C_dH^{-Am d_0d/((m-2)d_0+m(1-\alpha)-1)}$.
This bound also holds when $\delta_1(t)=R/2$.

Starting with \eqref{oa:baseline} and repeating this argument yields
radii bounded below by $C_nH^{-d_n}$, where
\begin{equation}\label{oa:improvement}
 d_{n+1}=\frac{Am d_0}{(m-2)d_0+m(1-\alpha)-1}\,d_n,
 \qquad n=0,1,2,\ldots.
\end{equation}
By \eqref{oa:conditionalindices} and \eqref{oa:contractionchoice},
\[
 Am-m+2=\frac{M-m}{M-1},\qquad
 0<\frac{Am d_0}{(m-2)d_0+m(1-\alpha)-1}<1.
\]
We may therefore fix a finite $n$ such that
\begin{equation}\label{oa:stop}
 Ad_n<1-\alpha-\frac1m.
\end{equation}
All constants after these finitely many steps are independent of $H,T$.
For each $t$, mass conservation provides a point
$\widetilde r(t)\in[0,R]$ with
$\rho(\widetilde r(t),t)=\bar\rho$.
Applying \eqref{oa:negativeholder} with $d=d_n$ at this point gives
\[
 \rho(0,t)^{\alpha-1+1/m}
 \le\bar\rho^{\alpha-1+1/m}+C_nH^{Ad_n}.
\]
Taking the supremum over $0\le t\le T$ and using the definition of $H$,
we obtain
\begin{equation}\label{oa:Habsorption}
 H^{1-\alpha-1/m}\le C_n\left(1+H^{Ad_n}\right).
\end{equation}
Together with \eqref{oa:stop}, this bounds $H$ independently of $T$.
Finally, \eqref{oa:onesided} gives
$\rho(r,t)\ge K^{-1}\rho(0,t)\ge C^{-1}$ for every $r\in[0,R]$.
This proves \eqref{oa:lowerestimate}.
\end{proof}

\subsection{Consequences for the effective velocity and the density}

With both density bounds established, the effective-velocity estimates
give the following bounds for $u$ and $\nabla\rho$.

\begin{proposition}\label{oa:bridge}
There is an exponent $m>2$ for which
\begin{align}
 &\sup_{0\le t\le T}
     \bigl(\norm{u(t)}_{L^m}+\norm{\nabla\rho(t)}_{L^m}\bigr)
 \notag\\
 &\qquad+\int_0^T
   \bigl(\norm{\nabla\rho}_{L^m}^m
        +\norm{\nabla\rho}_{L^2}^2
        +\norm{\nabla u}_{L^2}^2\bigr)\,dt\le C .    \label{oa:bridgeestimate}
\end{align}
\end{proposition}

\begin{proof}
The finite iteration above gives, for its fixed higher exponent
$m>1/(1-\alpha)$,
\[
 \sup_{0\le t\le T}\int \rho|V|^m\,dx\le C,
 \qquad
 \int_0^T\left(\int\rho F^m\,dx\right)^{1/m}\,dt\le C .
\]
Returning to \eqref{oa:Vmoment} and integrating in time, we obtain
\begin{equation}\label{oa:Vdamping}
 \int_0^T\int \rho^{\gamma-\alpha+1}|V|^m\,dx\,dt\le C .
\end{equation}
By Proposition \ref{oa:lower}, all powers of the density appearing
here are bounded above and below by positive constants.
Proposition \ref{oa:velocity} also gives
\[
 \sup_{0\le t\le T}\norm{u(t)}_{L^m}^m
 +\int_0^T\norm{u(t)}_{L^m}^m\,dt\le C ;
\]
for the time integral, use
$\int|u|^m\,dx\le C\int\rho^\alpha|u|^m/r^2\,dx
\le CD_{m-2}(t)$.
The relation \eqref{oa:Vdefinition} now implies
\[
 \sup_{0\le t\le T}\norm{\nabla\rho(t)}_{L^m}^m
 +\int_0^T\norm{\nabla\rho(t)}_{L^m}^m\,dt\le C.
\]
The remaining two time integrals in \eqref{oa:bridgeestimate}
follow from \eqref{eq:basic}, \eqref{eq:bd}, and the density bounds.

This completes the proof.
\end{proof}

Together, \eqref{oa:lowerestimate} and \eqref{oa:bridgeestimate} give
the bounds in \cite[(85)]{HMZ2025} with $2l=m>2$, uniformly for $T<T_*$.
The higher-order estimates and continuation argument of
\cite[Sections 3.6--3.7]{HMZ2025} then allow the local solution of
Lemma \ref{lem:local} to be extended past any finite $T_*$. Hence $T_*=\infty$.

\section{The gradient stress for \texorpdfstring{$\alpha=1$}{alpha = 1}}\label{ep:section}

In this section, we consider \eqref{eq:model} with $\sigma=1$ and
$\alpha=1$. Thus
\begin{equation}\label{ep:notation}
 V=U+2(\log\rho)_r,\qquad F=\frac{W^2}{r}.
\end{equation}
All the estimates below hold for every $\gamma>1$. We first establish
bounds on arbitrary finite time intervals, and then obtain the uniform
positive lower bound of the density. This order is needed for the
continuation argument.

\subsection{The upper bound of the density}

Fix an arbitrary interval $[0,T]$ on which the positive classical
solution from Lemma \ref{lem:local} exists.
In this section, a constant denoted by $C$ is independent of $T$, whereas
$C(T)$ may depend on $T$. We use $D_p$ from \eqref{eq:Dp} with
$\alpha=1$. By \eqref{eq:polar-gradient},
$|\nabla u|^2\ge |u|^2/r^2$.
We first estimate $\int\rho|u|^{p+2}\,dx$ and derive a uniform upper
bound for the density.

\begin{lemma}\label{ep:upper}
There is a constant $M>0$, independent of $T$, such that
\begin{equation}\label{ep:rhoupper}
 \rho(r,t)\le M\qquad (0\le r\le R,\ 0\le t\le T).
\end{equation}
Moreover, for every fixed $p\ge2$,
\begin{equation}\label{ep:velocitymoments}
 \sup_{0\le t\le T}\int \rho|u|^{p+2}\,dx
 +\int_0^T D_p(t)\,dt\le C_p,
\end{equation}
where $C_p$ depends on $p$ and the initial data, but not on $T$.
\end{lemma}

\begin{proof}
At $\alpha=1$, multiplying the second equation in \eqref{eq:model}
by $|u|^pu$, and using the first equation in \eqref{eq:model}, gives
\begin{align*}
 &\frac1{p+2}\frac{d}{dt}\int \rho|u|^{p+2}\,dx
   +2D_p
   +2p\int \rho|u|^{p-2}|(\nabla u)^Tu|^2\,dx\\
 &\qquad=\int \rho^\gamma
           \operatorname{div}(|u|^pu)\,dx
 \le D_p+C_p\int \rho^{2\gamma-1}|u|^p\,dx.
\end{align*}
For $p\ge2$, Young's inequality, \eqref{eq:polar-gradient}, and the
spherical estimate \eqref{eq:sphere} yield
\begin{align*}
 \int \rho^{2\gamma-1}|u|^p\,dx
 &=
 \int 
 \left(\frac{\rho|u|^{p+2}}{r^2}\right)^{(p-2)/p}
 \left(\frac{\rho|u|^2}{r^2}\right)^{2/p}
 \rho^{2\gamma-2}r^2\,dx\\
 &\le\varepsilon D_p+
 C_{\varepsilon,p}\int 
 \frac{\rho|u|^2}{r^2}\rho^{(\gamma-1)p}r^p\,dx
 \le\varepsilon D_p+C_{\varepsilon,p}D_0.
\end{align*}
For $p=2$, the first factor in the preceding integrand is understood
to be one, and the same estimate follows directly from
\eqref{eq:sphere}. Integrating in time and using \eqref{eq:basic}
proves \eqref{ep:velocitymoments}.

The radial effective velocity satisfies
\begin{equation}\label{ep:Veq}
 V_t+UV_r+\frac{\gamma}{2}\rho^{\gamma-1}V=\frac{\gamma}{2}\rho^{\gamma-1}U+F.
\end{equation}
Multiplying \eqref{ep:Veq} by $-\rho V_-^3$, integrating over
$\Omega$, and using the first equation in \eqref{eq:model}, gives
\[
 \frac14\frac{d}{dt}\int \rho V_-^4\,dx
   +\frac{\gamma}{4}\int \rho^\gamma V_-^4\,dx
 \le C\int \rho^\gamma|U|^4\,dx
 \le C\|\rho^{\gamma-1}r^2\|_{L^\infty}D_2.
\]
Here $V_-=\max\{-V,0\}$, and the term $F\ge0$ has a favorable
sign. It follows from \eqref{ep:velocitymoments} that
\begin{equation}\label{ep:negativefour}
 \sup_{0\le t\le T}\int \rho V_-^4\,dx\le C.
\end{equation}

By \eqref{eq:basic} and \eqref{eq:bd},
$\sqrt\rho$ is bounded in $H^1(\Omega)$, uniformly in time.
The one-dimensional trace inequality on $[R/2,R]$ therefore gives
$\rho(R,t)\le C$. Since $\rho_r=\rho(V-U)/2$, we obtain
\begin{align*}
 \rho(r,t)
 &\le \rho(R,t)+\frac12\int_0^R\rho(V_-+|U|)\,ds\\
 &\le C+
 C\left(\int_0^R\rho(V_-+|U|)^4s\,ds\right)^{1/4}
   \left(\int_0^R\rho s^{-1/3}\,ds\right)^{3/4}
 \le C.
\end{align*}
Indeed, the first integral is bounded by \eqref{ep:negativefour}
and \eqref{ep:velocitymoments}, whereas
$\rho s^{-1/3}\le Cs^{-2/3}$ follows from \eqref{eq:sphere}.
This proves \eqref{ep:rhoupper}.
\end{proof}

\subsection{Bounds for the negative radial velocities}

We next obtain lower bounds for the two radial velocities. The
nonnegative centrifugal term in \eqref{ep:Veq} is useful for this purpose.

\begin{lemma}\label{ep:negative}
There exist two constants $0<A<B$, independent of $T$, such that
\begin{equation}\label{ep:negativebounds}
 U(r,t)\ge-A,\qquad V(r,t)\ge-B
 \qquad (0\le r\le R,\ 0\le t\le T).
\end{equation}
\end{lemma}

\begin{proof}
The radial component of the second equation in \eqref{eq:model} can be written as
\begin{equation}\label{ep:Ueq}
 U_t+(2U-V)U_r
 =2U_{rr}+\frac{2U_r}{r}-\frac{2U}{r^2}
   +\frac{\gamma}{2}\rho^{\gamma-1}(U-V)+F.
\end{equation}
Choose $\kappa>0$ such that $e^{\kappa R}<2$, and set
$X=e^{\kappa r}U$. A direct calculation gives
\begin{align}
 &X_t+\left(2U-V+4\kappa-\frac2r\right)X_r-2X_{rr}\notag\\
 &\quad=2\kappa e^{-\kappa r}X^2+
 \left(2\kappa^2-\frac{2\kappa}{r}-\frac2{r^2}
       +\frac{\gamma}{2}\rho^{\gamma-1}\right)X\notag\\
 &\qquad-V\left(\kappa X+\frac{\gamma}{2}\rho^{\gamma-1}e^{\kappa r}\right)
       +e^{\kappa r}F.
 \label{ep:Xeq}
\end{align}
Fix $c\in(e^{\kappa R}/2,1)$. By taking $B$ sufficiently large and setting
$A=cB$, we can arrange that
\begin{gather}
 A>\|e^{\kappa r}U_{0,-}\|_{L^\infty},\qquad
 B>\|V_{0,-}\|_{L^\infty},\notag\\
 A>\frac{\gamma e^{\kappa R}M^{\gamma-1}}{2\kappa},\qquad
 A>\frac{e^{\kappa R}}2\left(B+\frac{2\kappa^2+\frac{\gamma}{2}M^{\gamma-1}}{\kappa}\right),\qquad
 A<B.
 \label{ep:AB}
\end{gather}
The choice is possible because $c>e^{\kappa R}/2$.

We claim that $X\ge-A$ and $V\ge-B$. Suppose first that $X$
reaches $-A$ from above at a first interior minimum, while
$V\ge-B$. At that point,
$X_t\le0$, $X_r=0$, $X_{rr}\ge0$, and
$\kappa X+\frac{\gamma}{2}\rho^{\gamma-1}e^{\kappa r}<0$.
Consequently, the right-hand side of \eqref{ep:Xeq} is bounded
from below by
\[
 2\kappa e^{-\kappa R}A^2-(2\kappa^2+\frac{\gamma}{2}M^{\gamma-1})A-\kappa BA>0.
\]
We have used the nonnegativity of the geometric terms and of $F$,
as well as
$-V(\kappa X+\frac{\gamma}{2}\rho^{\gamma-1}e^{\kappa r})\ge-\kappa BA$.
This contradicts the first downward contact with $-A$.

On the other hand, suppose that $V$ first reaches $-B$ at an
interior point $(r_0,t_0)$ while $X\ge-A$. At this point,
$V_t\le0$, $V_r=0$, and $U=e^{-\kappa r_0}X\ge-A$.
It follows from \eqref{ep:Veq} that
\[
V_t(r_0,t_0)
 =\frac{\gamma}{2}\rho(r_0,t_0)^{\gamma-1}
       \bigl(U(r_0,t_0)+B\bigr)+F(r_0,t_0)
 \ge\frac{\gamma}{2}\rho(r_0,t_0)^{\gamma-1}(B-A)>0,
\]
which is a contradiction. At the center,
$U=V=0$. At $r=R$, $U=W=0$ and
\begin{equation}\label{ep:Vboundary}
 V_t(R,t)+\frac{\gamma}{2}\rho(R,t)^{\gamma-1}V(R,t)=0.
\end{equation}
Thus neither boundary can be the first point where either lower
bound is lost. This proves the claim and hence
\eqref{ep:negativebounds}.
\end{proof}

\subsection{Finite-time bounds for the physical velocity}

Using these lower bounds, we obtain upper bounds for the physical
velocity on every finite time interval by a first-contact argument.

\begin{lemma}\label{ep:physical}
There exists a constant $C(T)>0$, depending only on $T$, $\gamma$,
$R$ and the initial data, such that
\begin{equation}\label{ep:physicalbounds}
 \sup_{0\le t\le T}\|u(t)\|_{L^\infty}\le C(T).
\end{equation}
\end{lemma}

\begin{proof}
At $\alpha=1$, equation \eqref{eq:angular} can be written as
\begin{equation}\label{ep:Weq}
 W_t+\left(2U-V-\frac2r\right)W_r-2W_{rr}
 =\left(-\frac Ur-\frac2{r^2}\right)W.
\end{equation}
By Lemma \ref{ep:negative},
\[
 -\frac Ur-\frac2{r^2}
 \le\frac Ar-\frac2{r^2}\le\frac{A^2}{8}.
\]
Choose
\[
 C_{w,0}=
 \max\{1,\|W_0\|_{L^\infty},R\|W_{0,r}\|_{L^\infty}\}.
\]
For $\varepsilon>0$, suppose that $W$ first touches
$(C_{w,0}+\varepsilon)e^{(A^2/8+\varepsilon)t}$ from below.
Since the initial inequality is strict and $W(0,t)=W(R,t)=0$,
the contact point lies in $0<r<R$ at a positive time.
At this point, $W>0$, $W_r=0$, $W_{rr}\le0$, and
\eqref{ep:Weq} gives
\[
 \left(\frac{A^2}{8}+\varepsilon\right)W
 \le W_t
 \le\left(-\frac Ur-\frac2{r^2}\right)W
 \le\frac{A^2}{8}W,
\]
a contradiction. Applying the same argument to $-W$ and letting
$\varepsilon\to0$, we obtain
\begin{equation}\label{ep:Cw}
 |W(r,t)|\le C_w(t):=C_{w,0}e^{A^2t/8}.
\end{equation}

At an interior positive maximum of $U$, equation
\eqref{ep:Ueq} implies
\[
 U_t\le-\frac{2U}{r^2}+\frac{\gamma}{2}M^{\gamma-1}(U+B)+\frac{C_w(t)^2}{r}.
\]
For $U\ge1$,
\[
 -\frac{2U}{r^2}+\frac{C_w(t)^2}{r}
 \le\max_{z\ge0}\{-2Uz^2+C_w(t)^2z\}
 =\frac{C_w(t)^4}{8U}\le\frac{C_w(t)^4}{8}.
\]
We may therefore take $L$ to be the solution of
\begin{equation}\label{ep:L}
 L'=\frac{\gamma}{2}M^{\gamma-1}(L+B)+\frac{C_w(t)^4}{8},\qquad
 L(0)=\max\{1,\|U_{0,+}\|_{L^\infty}\}.
\end{equation}
This linear equation defines $L$ on $[0,\infty)$, and $L$ is
finite on every bounded time interval.
To prove $U\le L$, let $L_\varepsilon$ solve the same equation
with an additional term $+\varepsilon$ on the right-hand side
and initial value $L(0)+\varepsilon$. Since $U(0,t)=U(R,t)=0$
and $L_\varepsilon\ge1$, a first contact of $U$ with
$L_\varepsilon$ from below must occur at an interior point.
There, $U=L_\varepsilon$, $U_r=0$, $U_{rr}\le0$, and the
preceding estimates give
\[
 L_\varepsilon'\le U_t
 \le\frac{\gamma}{2}M^{\gamma-1}(L_\varepsilon+B)
       +\frac{C_w(t)^4}{8}
 =L_\varepsilon'-\varepsilon,
\]
which is impossible. Letting $\varepsilon\to0$ gives $U_+\le L$.
Together with $U\ge-A$, this proves \eqref{ep:physicalbounds}
with $C(T)=A+C_w(T)+L(T)$.
\end{proof}

\subsection{The angular barrier and the effective velocity}

The following estimate gives a linear bound for the centrifugal
force in terms of the positive effective velocity.

\begin{lemma}\label{ep:localbarrier}
Let
\begin{equation}\label{ep:K}
 K(t)=1+\sup_{0\le s\le t}\|V_+(s)\|_{L^\infty}.
\end{equation}
Then
\begin{equation}\label{ep:linearforce}
 \left\|\frac{W(t)^2}{r}\right\|_{L^\infty}
 \le C_w(t)^2\left(\frac2R+4A+1+K(t)\right).
\end{equation}
Here $C_w$ is the function in \eqref{ep:Cw}.
\end{lemma}

\begin{proof}
Fix $t\in[0,T]$ and choose
\begin{equation}\label{ep:delta}
 \delta=\min\left\{\frac R2,\frac1{K(t)+4A+1}\right\}.
\end{equation}
The number $\delta$ is held fixed on the whole cylinder
$0<r<\delta$, $0\le s\le t$. On this cylinder, put
$Z(r,s)=W(r,s)/\sqrt r$. Equation \eqref{ep:Weq} gives
\begin{equation}\label{ep:Zeq}
 Z_s+\left(2U-V-\frac4r\right)Z_r-2Z_{rr}
 =\left(\frac{V/2-2U}{r}-\frac3{2r^2}\right)Z.
\end{equation}
Since $V(r,s)\le K(t)$ and $U(r,s)\ge-A$, its zeroth-order
coefficient satisfies
\begin{equation}\label{ep:reaction}
 \frac{V/2-2U}{r}-\frac3{2r^2}
 \le\frac{(K(t)/2+2A)r-3/2}{r^2}
 \le-\frac1{r^2}<0.
\end{equation}
At $r=\delta$,
\[
 |Z(\delta,s)|\le C_w(s)\delta^{-1/2}
 \le C_w(t)\delta^{-1/2}.
\]
Since $W_0(0)=0$, $\delta\le R$ and
$C_{w,0}\ge R\|W_{0,r}\|_{L^\infty}$, we also have
\[
 |Z(r,0)|
 \le\|W_{0,r}\|_{L^\infty}\sqrt r
 \le\frac{R\|W_{0,r}\|_{L^\infty}}{\sqrt\delta}
 \le\frac{C_{w,0}}{\sqrt\delta}
 \le\frac{C_w(t)}{\sqrt\delta}.
\]
Moreover, $Z$ can extend
continuously to $r=0$, with $Z(0,s)=0$.

For $\varepsilon>0$, suppose that $Z$ first touches the constant
$C_w(t)\delta^{-1/2}+\varepsilon$ from below. By the initial
and boundary estimates above, the contact point lies in
$0<r<\delta$ at a positive time. There, $Z>0$, $Z_s\ge0$,
$Z_r=0$ and $Z_{rr}\le0$. Equations \eqref{ep:Zeq} and
\eqref{ep:reaction} therefore give
\[
 0\le Z_s-2Z_{rr}
 =\left(\frac{V/2-2U}{r}-\frac3{2r^2}\right)Z
 \le-\frac Z{r^2}<0,
\]
a contradiction. Applying the same argument to $-Z$ and letting
$\varepsilon\to0$, we obtain
\[
 |W(r,s)|\le C_w(t)\sqrt{r/\delta}
 \qquad(0\le r\le\delta,\ 0\le s\le t).
\]
Thus, at time $t$, $W^2/r\le C_w(t)^2/\delta$ in the inner
region. The same estimate in $r\ge\delta$ follows from
\eqref{ep:Cw}. Finally,
\[
 \delta^{-1}
 =\max\{2/R,K(t)+4A+1\}
 \le 2/R+K(t)+4A+1,
\]
which proves \eqref{ep:linearforce}.
\end{proof}

The linear bound \eqref{ep:linearforce} closes the estimate for the
positive effective velocity by Gronwall's inequality. It also yields
a positive lower bound for the density on every finite time interval.

\begin{proposition}\label{ep:finitetime}
We have
\begin{equation}\label{ep:finitebounds}
 \sup_{0\le t\le T}
 \left(\|u(t)\|_{L^\infty}+\|V(t)\|_{L^\infty}
 +\left\|\frac{W(t)^2}{r}\right\|_{L^\infty}
 +\|\nabla\log\rho(t)\|_{L^\infty}\right)
 \le C(T).
\end{equation}
Moreover,
\begin{equation}\label{ep:finitelower}
 \rho(r,t)\ge\bar\rho\,e^{-RC(T)}>0
 \qquad(0\le t\le T),
 \quad
 \bar\rho=\frac1{|\Omega|}\int \rho_0\,dx.
\end{equation}
Consequently, the solution extends globally in time. 
\end{proposition}

\begin{proof}
Along material trajectories, \eqref{ep:Veq} implies
\[
 D_tV_+\le \frac{\gamma}{2}M^{\gamma-1}U_++F.
\]
Using Lemma \ref{ep:localbarrier} at each time $s$, and taking
the running spatial maximum, we obtain
\begin{align}
 K(t)\le{}&1+\|V_{0,+}\|_{L^\infty}\notag\\
 &+\int_0^t\left[
 \frac{\gamma}{2}M^{\gamma-1}L(s)+C_w(s)^2
       \left(\frac2R+4A+1+K(s)\right)\right]\,ds.
 \label{ep:gronwall}
\end{align}
Here $L$ is the solution of \eqref{ep:L}, constructed in the proof
of Lemma \ref{ep:physical}.
Both $L$ and $C_w$ are bounded on every finite interval.
Gronwall's inequality therefore yields $K(t)\le C(T)$ for
$0\le t\le T$. Together with Lemmas \ref{ep:negative} and
\ref{ep:physical}, and \eqref{ep:linearforce}, this proves all
the velocity bounds in \eqref{ep:finitebounds}.
The identity $(\log\rho)_r=(V-U)/2$ gives the last bound.

For every $t$, mass conservation and continuity provide a point
$\widetilde r(t)\in[0,R]$ such that
$\rho(\widetilde r(t),t)=\bar\rho$. Hence
\[
 |\log\rho(r,t)-\log\bar\rho|
 \le\int_0^R|(\log\rho)_s|\,ds\le RC(T),
\]
which proves \eqref{ep:finitelower}.

If $T_*<\infty$, the preceding estimates bound $\rho$, $\rho^{-1}$,
$u$ and $\nabla\rho$ on $[0,T_*)$. Together with \eqref{eq:basic},
they imply the bound in \cite[(85)]{HMZ2025} with $2l=3$, uniformly for $T<T_*$.
The continuation argument following Proposition \ref{oa:bridge}
therefore gives global existence.
\end{proof}

\subsection{The time-uniform lower bound of the density}

We now work with the global solution obtained above.
Only the upper density bound is needed in applying
Lemma \ref{lem:decay}; in particular, it gives
\begin{equation}\label{ep:decay}
 \int \rho|u|^2\,dx+\|\rho(t)-\bar\rho\|_{L^2}^2
 +\int_t^\infty D_0(s)\,ds\le Ce^{-\lambda t}.
\end{equation}
We use weighted estimates for $W$ to control the centrifugal term
in the estimate for $\rho^{1/3}V$ in $L^3$.

\begin{lemma}\label{ep:angularmoments}
The centrifugal force satisfies
\begin{equation}\label{ep:Ftime}
 \int_0^\infty
 \left(\int\rho F(t)^3\,dx\right)^{1/3}\,dt\le C.
\end{equation}
\end{lemma}

\begin{proof}
Multiplying \eqref{eq:angular} with $\alpha=1$ by
$|W|^{q-2}Wr^{-k}$, where $q\ge2$ and $0<k<q$, and using the first
equation in \eqref{eq:model}, we obtain
\begin{align}
 &\frac1q\frac{d}{dt}\int \rho|W|^qr^{-k}\,dx\notag\\
 &\quad+
 2\int \rho|W|^{q-2}r^{-k}
 \left[(q-1)W_r^2-kW_r\frac Wr+\frac{W^2}{r^2}\right]\,dx
 \notag\\
 &\qquad=
 -\left(1+\frac kq\right)
       \int \rho U|W|^qr^{-k-1}\,dx.
 \label{ep:angularidentity}
\end{align}
The quadratic form is strictly positive definite whenever
$k^2<4(q-1)$.

Take first $q=2$ and $k=3/2$. By Young's inequality,
\begin{align*}
 \int \rho|U|W^2r^{-5/2}\,dx
 &\le\varepsilon\int \rho W^2r^{-7/2}\,dx
   +C_\varepsilon\int \rho U^2W^2r^{-3/2}\,dx\\
 &\le\varepsilon\int \rho W^2r^{-7/2}\,dx
   +C_\varepsilon D_2.
\end{align*}
In the last step we used $r^{1/2}\le R^{1/2}$ and
\eqref{eq:polar-gradient}. The quadratic form in
\eqref{ep:angularidentity} absorbs the first term. Therefore
\begin{equation}\label{ep:weightedtwo}
 \sup_{t\ge0}\int \rho W^2r^{-3/2}\,dx\le C.
\end{equation}
Writing $J_j(t)=\int \rho F^j\,dx$, we infer from
H\"older's inequality and \eqref{ep:decay} that
\begin{equation}\label{ep:Jone}
 J_1
 \le\left(\int \rho W^2\,dx\right)^{1/3}
      \left(\int \rho W^2r^{-3/2}\,dx\right)^{2/3}
 \le Ce^{-\lambda t/3}.
\end{equation}

Next, take $q=8$ and $k=2$ in \eqref{ep:angularidentity}.
By the positive definiteness of the quadratic form, we have
\[
2\int \rho|W|^6r^{-2}
\left(7W_r^2-2W_r\frac Wr+\frac{W^2}{r^2}\right)\,dx
\ge c\int \rho|W|^8r^{-4}\,dx=cJ_4(t).
\]
For the transport term on the right-hand side of \eqref{ep:angularidentity},
we use the estimate
\begin{align*}
 \int \rho|U||W|^8r^{-3}\,dx
 &\le J_4^{1/2}
 \left(\int \rho U^2W^8r^{-2}\,dx\right)^{1/2}\\
 &\le J_4^{1/2}
 \left(\int \rho|u|^{10}r^{-2}\,dx\right)^{1/2}
 \le\varepsilon J_4+C_\varepsilon D_8.
\end{align*}
Using \eqref{ep:velocitymoments}, we obtain
\begin{equation}\label{ep:Jfour}
 \int_0^\infty J_4(t)\,dt\le C.
\end{equation}
H\"older's inequality gives
\[
 \left(\int\rho F^3\,dx\right)^{1/3}
 \le J_1^{1/9}J_4^{2/9}.
\]
Consequently, by \eqref{ep:Jone}, \eqref{ep:Jfour}, and
H\"older's inequality in time,
\[
 \int_0^\infty\left(\int\rho F^3\,dx\right)^{1/3}\,dt
 \le
 \left(\int_0^\infty J_1^{1/7}\,dt\right)^{7/9}
 \left(\int_0^\infty J_4\,dt\right)^{2/9}\le C.
\]
This proves \eqref{ep:Ftime}.
\end{proof}

We now use \eqref{ep:Ftime} to bound $\int\rho|V|^3\,dx$ uniformly
in time and deduce a positive lower bound for the density independent
of time.

\begin{lemma}\label{ep:uniformlower}
There is a constant $C>0$, independent of time, such that
\begin{equation}\label{ep:uniformdensity}
 C^{-1}\le\rho(r,t)\le C\qquad (0\le r\le R,\ t\ge0).
\end{equation}
In addition,
\begin{equation}\label{ep:Vthree}
 \sup_{t\ge0}\int \rho|V|^3\,dx
 +\int_0^\infty\int \rho^\gamma|V|^3\,dx\,dt\le C,
\end{equation}
and
\begin{equation}\label{ep:densityconvergence}
 \sup_{t\ge0}\|\nabla\rho(t)\|_{L^3}\le C,\qquad
 \|\rho(t)-\bar\rho\|_{L^\infty}\le Ce^{-\lambda t/8}.
\end{equation}
\end{lemma}

\begin{proof}
Multiplying \eqref{ep:Veq} by $\rho|V|V$, integrating, and using
the first equation in \eqref{eq:model} and Young's inequality, we obtain
\begin{equation}\label{ep:Yineq}
\begin{split}
 &\frac13\frac{d}{dt}\int\rho|V|^3\,dx
       +c\int\rho^\gamma|V|^3\,dx\\
 &\quad\le C(D_0+D_2)
       +\left(\int\rho F^3\,dx\right)^{1/3}
        \left(\int\rho|V|^3\,dx\right)^{2/3}.
\end{split}
\end{equation}
Here we used the upper density bound and \eqref{eq:polar-gradient}
to estimate the pressure remainder:
\[
 \int\rho^\gamma|U|^3\,dx
 \le C\int\rho(|u|^2+|u|^4)\,dx
 \le C(D_0+D_2).
\]
Multiplying \eqref{ep:Yineq} by
$\left(1+\int\rho|V|^3\,dx\right)^{-2/3}$ and discarding the
nonnegative damping term gives
\[
 \frac{d}{dt}\left(1+\int\rho|V|^3\,dx\right)^{1/3}
 \le C(D_0+D_2)+\left(\int\rho F^3\,dx\right)^{1/3}.
\]
Integration in time, together with \eqref{eq:basic},
\eqref{ep:velocitymoments} and \eqref{ep:Ftime}, yields
$\sup_{t\ge0}\int\rho|V|^3\,dx\le C$.
Integrating \eqref{ep:Yineq} again and retaining the damping term
then proves \eqref{ep:Vthree}.

Since
\[
 \partial_r\rho^{1/3}=\frac16\rho^{1/3}(V-U),
\]
\eqref{ep:Vthree}, \eqref{eq:basic}, and
\eqref{ep:velocitymoments} imply
$\sup_{t\ge0}\|\nabla\rho^{1/3}(t)\|_{L^3}\le C$.
The density upper bound therefore gives
$\sup_{t\ge0}\|\nabla\rho(t)\|_{L^3}\le C$.
The two-dimensional interpolation inequality and
\eqref{ep:decay} yield
\begin{align*}
 \|\rho(t)-\bar\rho\|_{L^\infty}
 &\le C\|\rho(t)-\bar\rho\|_{L^2}^{1/4}
       \|\nabla\rho(t)\|_{L^3}^{3/4}
       +C\|\rho(t)-\bar\rho\|_{L^2}\\
 &\le Ce^{-\lambda t/8}.
\end{align*}
This proves \eqref{ep:densityconvergence}.

The solution is already global by Proposition
\ref{ep:finitetime}. Hence we may choose a fixed $T_0$,
depending only on the data, such that
$\rho(r,t)\ge\bar\rho/2$ for $t\ge T_0$.
On $[0,T_0]$, the positive lower bound follows from
\eqref{ep:finitelower}. Combining these two estimates proves
\eqref{ep:uniformdensity}.
\end{proof}

For later use, the endpoint estimates also give
\begin{equation}\label{ep:uniforminputs}
 \begin{split}
 &\sup_{t\ge0}
 \left(\|u(t)\|_{L^3}+
       \|\nabla\rho(t)\|_{L^2}+
       \|\nabla\rho(t)\|_{L^3}\right)\\
 &\quad+
 \int_0^\infty
 \left(\|\nabla u(t)\|_{L^2}^2+
       \|\nabla\rho(t)\|_{L^2}^2+
       \|\nabla\rho(t)\|_{L^3}^3\right)\,dt\le C.
 \end{split}
\end{equation}
Indeed, the quadratic bounds follow from the basic energy and
BD entropy after using \eqref{ep:uniformdensity}. Moreover,
\[
 \int |\nabla\rho|^3\,dx
 \le C\int \rho^\gamma|V|^3\,dx
      +C\int \rho|u|^3\,dx
 \le C\int \rho^\gamma|V|^3\,dx+C(D_0+D_2).
\]
Integration in time and \eqref{ep:Vthree} prove the remaining
bound.

\section{The symmetric stress}\label{sy:section}

In this section, we take $\sigma=0$ and $1/2<\alpha\le1$.
We first derive estimates on an arbitrary interval $[0,T]$ on which
the classical solution exists and the density is positive. With $F$ as
in \eqref{oa:F}, equation \eqref{eq:radial-effective} becomes
\[
 V_t+UV_r=\frac{\gamma}{2\alpha}\rho^{\gamma-\alpha}(U-V)+F.
\]
Since the radial component of the divergence of
$\rho^\alpha(\nabla u-\nabla u^\top)$ vanishes,
the radial momentum equation can be written as
\begin{equation}\label{sy:Ueq}
\begin{split}
 U_t+\bigl((1+\alpha)U-\alpha V\bigr)U_r
 &=2\alpha\rho^{\alpha-1}\left(U_{rr}+\frac{U_r}{r}-\frac{U}{r^2}\right)\\
 &\quad+\left((\alpha-1)\frac Ur
       -\frac{\gamma}{2\alpha}\rho^{\gamma-\alpha}\right)(V-U)+F.
\end{split}
\end{equation}
The tangential component of the momentum equation is now
\begin{equation}\label{sy:angular-stress}
 \rho\left(W_t+UW_r+\frac{UW}{r}\right)
 =\frac1{r^2}\partial_r\left[r^2\rho^\alpha
                    \left(W_r-\frac Wr\right)\right].
\end{equation}
Expanding the right-hand side and using
$2(\rho^\alpha)_r=\rho(V-U)$, we also have
\begin{equation}\label{sy:angular}
 \rho\left(W_t+UW_r+\frac{U+V}{2r}W\right)
 =\frac1r(r\rho^\alpha W_r)_r-\rho^\alpha\frac{W}{r^2}.
\end{equation}

\subsection{Lower bounds for the radial velocities}

We first estimate the negative parts of $U$ and $V$ without assuming
an upper bound for the density.

\begin{lemma}\label{sy:negative}
There exists a constant $A>0$, independent of $T$, such that
\[
 U(r,t)\ge-A,\qquad V(r,t)\ge-A
 \qquad(0\le r\le R,\ 0\le t\le T).
\]
\end{lemma}

\begin{proof}
Fix $0<\delta<\alpha$ and put $f(r)=\kappa r^{1-\delta}$.
We take $\kappa>0$ sufficiently small that
\[
 rf'(r)\le1,\qquad
 (1+\alpha)e^{-f(R)}-\alpha-\frac{1-\alpha}{1-\delta}>0 .
\]
Indeed, $rf'(r)\le\kappa(1-\delta)R^{1-\delta}$, and the expression in
the second inequality tends to $(\alpha-\delta)/(1-\delta)>0$ as $\kappa\to0$.
By \eqref{eq:sphere}, $r^\delta\rho^{\gamma-\alpha}\le C$ uniformly
in time. Since $f'(r)=\kappa(1-\delta)r^{-\delta}$ and $f(r)\le f(R)$,
we may choose $A>0$, independently of $T$, so that
\[
 A\alpha f'(r)>\frac{\gamma}{2\alpha}\rho^{\gamma-\alpha}e^{f(r)}
 \qquad (0<r\le R,\ 0\le t\le T).
\]

Set $X=e^fU$ and increase $A$ if necessary so that $X,V>-A$
initially. Suppose there is a first time when either function reaches
$-A$. At that time both are at least $-A$ everywhere.
If $X=-A$ at an interior point, then $X_t\le0$, $X_r=0$ and
$X_{rr}\ge0$ there. At this point $U_r=-e^{-f}f'X$ and
$U_{rr}=e^{-f}\bigl(X_{rr}+(f'^2-f'')X\bigr)$, so \eqref{sy:Ueq} gives
\begin{align*}
X_t-2\alpha\rho^{\alpha-1}X_{rr}
&=2\alpha\rho^{\alpha-1}\left(f'^2-f''-\frac{f'}r-\frac1{r^2}\right)X
   +\frac{\gamma}{2\alpha}\rho^{\gamma-\alpha}X+e^fF\\
&\quad+\left((1+\alpha)f'+\frac{1-\alpha}{r}\right)e^{-f}X^2\\
&\quad+\left(-X\left(\alpha f'+\frac{1-\alpha}{r}\right)
          -\frac{\gamma}{2\alpha}\rho^{\gamma-\alpha}e^f\right)V .
\end{align*}
Since $rf'\le1$ and $f''+f'/r=(1-\delta)f'/r\ge0$, the factor
multiplying $X=-A$ in the first term is nonpositive. The coefficient
of $V$ is positive by the choice of $A$ and $1-\alpha\ge0$.
We may therefore replace $V$ by $-A$ to obtain a lower bound.
The two pressure terms then sum to
$\frac{\gamma}{2\alpha}\rho^{\gamma-\alpha}A(e^f-1)\ge0$.
Dropping this contribution, the first term, and $e^fF\ge0$ gives
\begin{align*}
 X_t-2\alpha\rho^{\alpha-1}X_{rr}
 &\ge A^2\left[\bigl((1+\alpha)e^{-f}-\alpha\bigr)f'
              -(1-\alpha)\frac{1-e^{-f}}{r}\right]\\*
 &\ge A^2 f'\left((1+\alpha)e^{-f(R)}
                -\alpha-\frac{1-\alpha}{1-\delta}\right)>0.
\end{align*}
Here we used $1-e^{-f}\le f=rf'/(1-\delta)$ and the choice of $\kappa$.
This contradicts $X_t\le0$ and $X_{rr}\ge0$ at the contact point.

If $V=-A$ at an interior point, then $U=e^{-f}X\ge-Ae^{-f}>-A$. At first contact $V_t\le0$ and $V_r=0$, whereas
\eqref{eq:radial-effective} yields
$V_t=\frac{\gamma}{2\alpha}\rho^{\gamma-\alpha}(U+A)+F>0$, again a
contradiction. At the center, regularity and axisymmetry give $U(0,t)=V(0,t)=0$.
At $r=R$, the boundary condition gives $U=W=0$, hence $X=F=0$.
If $V(R,t)=-A$ at first contact, \eqref{eq:radial-effective} gives
$V_t(R,t)=\frac{\gamma}{2\alpha}\rho(R,t)^{\gamma-\alpha}A>0$,
again contradicting $V_t(R,t)\le0$. Thus $X,V\ge-A$ throughout
$[0,R]\times[0,T]$, and $U=e^{-f}X\ge-Ae^{-f}\ge-A$.
\end{proof}

\subsection{An upper bound for the density}

\begin{lemma}\label{sy:upper}
There exists a constant $C>0$, independent of $T$, such that
\[
 \rho(r,t)\le C\qquad(0\le r\le R,\ 0\le t\le T).
\]
\end{lemma}

\begin{proof}
Set $M(r,t)=\int_0^r s\rho(s,t)\dd s$. Integrating the continuity
equation and using \eqref{eq:effective} with $M_r=r\rho$ gives
\[
 M_t=-r\rho U
 =2\alpha\rho^{\alpha-1}\left(M_{rr}-\frac{M_r}{r}\right)-VM_r.
\]
With $A$ from Lemma \ref{sy:negative}, set $\Psi(r)=\int_0^r s w(s)\dd s$,
where $2\alpha w^{\alpha-2}w'=-A$ and $w(0)>0$. Thus
\[
 w(r)=
 \begin{cases}
 \left(w(0)^{\alpha-1}
        +\dfrac{A(1-\alpha)}{2\alpha}r\right)^{-1/(1-\alpha)},
       &\alpha<1,\\[2pt]
 w(0)e^{-Ar/2},&\alpha=1.
 \end{cases}
\]
As $w(0)\to\infty$, we have $\Psi(R)\to\infty$. Indeed, for $\alpha<1$
the integrand $s w(s)$ increases to a multiple of $s^{1-1/(1-\alpha)}$,
whose integral diverges at zero since $\alpha>1/2$. For $\alpha=1$,
$\Psi(R)$ is proportional to $w(0)$. Choose $w(0)$ so that
$\Psi(R)/R^2\ge\norm{\rho_0}_\infty/2$. Since $w$ is decreasing,
$\Psi(r)/r^2\ge\Psi(R)/R^2$, and hence
\[
 M(r,0)\le\tfrac12\norm{\rho_0}_{L^\infty}r^2\le\Psi(r)
 \qquad(0\le r\le R).
\]

For $\varepsilon>0$, the comparison $M<\Psi+\varepsilon e^t$ is strict
initially and at both endpoints, since $M(0,t)=0$ and $M(R,t)=M(R,0)$. At a first interior contact,
$M_r=\Psi'$ gives $\rho=w$, and $M_{rr}\le\Psi''$. Since $V\ge-A$,
\[
 0<\varepsilon e^t\le M_t
 \le 2\alpha w^{\alpha-1}\left(\Psi''-\frac{\Psi'}{r}\right)
     +A\Psi'
 =r\bigl(2\alpha w^{\alpha-1}w'+Aw\bigr)=0.
\]
Thus $M(r,t)\le\Psi(r)\le w(0)r^2/2$ upon letting $\varepsilon\to0$.

For $0<r\le R$, the weighted mean value theorem gives $s\in[r/2,r]$ with
\[
 \rho(s,t)=\frac{M(r,t)-M(r/2,t)}{\int_{r/2}^r \xi\dd\xi}
 \le\frac{8M(r,t)}{3r^2}\le\frac43 w(0).
\]
Integrating from $s$ to $r$ and applying Cauchy--Schwarz with
$\int_s^r\dd\xi/\xi\le\log 2$, we obtain
\[
 \rho(r,t)^{\alpha-1/2}
 \le \rho(s,t)^{\alpha-1/2}
    +\left(\frac{\log 2}{2\pi}\right)^{1/2}
       \norm{\nabla\rho^{\alpha-1/2}}_2 \le C
\]
by \eqref{eq:sphere}. Since $\alpha>1/2$, this bounds $\rho$ uniformly
for $r>0$, and continuity gives the same bound at $r=0$.
\end{proof}

\subsection{Finite-time positivity and continuation}

\begin{proposition}\label{sy:finitetime}
We have
\[
 \sup_{0\le t\le T}
 \bigl(\norm{u(t)}_{L^\infty}+\norm{V(t)}_{L^\infty}
       +\norm{F(t)}_{L^\infty}+\norm{\nabla\rho(t)}_{L^\infty}
       +\norm{\rho(t)^{-1}}_{L^\infty}\bigr)\le C(T).
\]
In particular, the classical solution exists globally in time.
\end{proposition}

\begin{proof}
By Lemma \ref{sy:upper}, $\rho^{\alpha-1}\ge c>0$ and
$\rho^{\gamma-\alpha}\le C$.
We first put $Z=W/\sqrt r$. From \eqref{sy:angular}, we obtain
\[
 Z_t+\left(\frac{3U-V}{2}-\frac{2\rho^{\alpha-1}}{r}\right)Z_r-\rho^{\alpha-1}Z_{rr}
 =\left(-\frac{3\rho^{\alpha-1}}{4r^2}-\frac{5U+V}{4r}\right)Z .
\]
Lemma \ref{sy:negative} implies that the coefficient of $Z$ on
the right is bounded above by $-3c/(4r^2)+3A/(2r)$, and hence by a
constant $C\ge0$. Since $W=O(r)$ at the center, $Z$ extends
continuously to zero at both endpoints. At a first upper contact
with $(\norm{Z(\cdot,0)}_{L^\infty}+\varepsilon)e^{(C+\varepsilon)t}$,
the drift term vanishes and the diffusion has nonpositive sign.
The equation would give $(C+\varepsilon)Z\le CZ$, which is
impossible. Applying the same comparison to $-Z$ and letting
$\varepsilon\to0$, we find $\norm{F(t)}_{L^\infty}
=\norm{Z(t)}_{L^\infty}^2\le C(T)$.

Next, set $b=(1-\alpha)/\alpha\in[0,1)$ and $X=r^{-b}U$.
By \eqref{sy:Ueq},
\begin{align*}
&\quad X_t+\left((1+\alpha)U-\alpha V
          -\frac{2\alpha(2b+1)\rho^{\alpha-1}}{r}\right)X_r\\
&=2\alpha\rho^{\alpha-1}X_{rr}
  -\left(\frac{2\alpha(1-b^2)\rho^{\alpha-1}}{r^2}+\frac{bU}{r}\right)X\\
&\quad+\frac{\gamma}{2\alpha}\rho^{\gamma-\alpha}(X-r^{-b}V)+r^{-b}F .
\end{align*}
The $VX/r$ term cancels because $\alpha b=1-\alpha$. Since $U=O(r)$
at the center and $b<1$, $X$ extends continuously to zero there,
and $X(R,t)=0$. At a spatial maximum with $X\ge1$, we have
$X_r=0$, $X_{rr}\le0$, and $U=r^bX>0$. Thus $-bUX/r\le0$.
Using $V\ge-A$ and $F\le C(T)$ in the remaining terms gives
\[
 X_t\le\left(C-\frac{2\alpha(1-b^2)c}{r^2}\right)X+C(T)r^{-b}
 \le K_T X .
\]
Here $K_T\ge0$ is finite. A first contact with
$(1+\norm{X(\cdot,0)}_\infty)e^{(K_T+1)t}$ would occur at an interior
maximum with $X_t\ge(K_T+1)X$, a contradiction. Hence $X_+$ is
bounded on $[0,T]$, and $U_+=r^bX_+\le R^bX_+$.

To bound $V_+$, let $r(t)$ be the radial coordinate of a fluid particle,
so that $\dot r(t)=U(r(t),t)$.
Where $V>0$, the term $-\frac{\gamma}{2\alpha}\rho^{\gamma-\alpha}V$
in \eqref{eq:radial-effective} is nonpositive. Thus along the trajectory,
for almost every $t$,
\[
 \frac{\dd}{\dd t}V_+(r(t),t)
 \le \frac{\gamma}{2\alpha}\rho^{\gamma-\alpha}U_+ + F \le C(T).
\]
These radial coordinates cover $[0,R]$ since $U$ vanishes at both endpoints.
Integrating gives $\norm{V_+(t)}_\infty\le\norm{V_+(\cdot,0)}_\infty+tC(T)$.

If $\alpha<1$, we have
$(\rho^{\alpha-1})_r=(\alpha-1)(V-U)/(2\alpha)$;
if $\alpha=1$, $(\log\rho)_r=(V-U)/2$.
Integrating from a point where the density equals $\bar\rho$ gives
$\rho^{-1}\le C(T)$. Also,
$\rho_r=\rho^{2-\alpha}(V-U)/(2\alpha)$, so $\norm{\nabla\rho}_{L^\infty}
\le C(T)$. All the constants in these estimates can be chosen
bounded when $T$ ranges over a bounded interval.

If $T_*<\infty$, these bounds and \eqref{eq:basic} give the estimate
in \cite[(85)]{HMZ2025} with $2l=3$, uniformly for $T<T_*$. Using
Korn's inequality and standard elliptic estimates for the Lam\'e system,
we conclude as in Proposition \ref{ep:finitetime} that $T_*=\infty$.
\end{proof}

\subsection{Weighted estimates for the tangential velocity}

For the global solution obtained in Proposition \ref{sy:finitetime},
the energy identity \eqref{eq:energy-identity} shows that the relative
energy \eqref{eq:relative-energy} is nonincreasing. By Lemma \ref{sy:upper} and
\eqref{eq:relative-equivalence}, we have
\[
 C^{-1}\left(\int\rho|u|^2\dd x+\norm{\rho-\bar\rho}_2^2\right)
 \le\mathcal E(t)\le C\left(\int\rho|u|^2\dd x+\norm{\rho-\bar\rho}_2^2\right).
\]
The first term is bounded by $CD_0(t)$ by \eqref{eq:polar-gradient}.
To estimate the second term, set
$s=(\gamma+\alpha-1)/2$ and $\rho_*(t)=\bigl(\overline{\rho^s}(t)\bigr)^{1/s}$.
By Jensen's inequality for $s\ge1$ and
$\rho^s\ge\norm{\rho}_\infty^{s-1}\rho$ for $s<1$, mass conservation
and the upper density bound give $C^{-1}\le\rho_*(t)\le C$.
Therefore $|z-\rho_*(t)|\le C|z^s-\rho_*(t)^s|$ for
$0\le z\le\norm{\rho}_\infty$.
Poincar\'e's inequality gives
\[
 \norm{\rho-\bar\rho}_2\le\norm{\rho-\rho_*(t)}_2
 \le C\norm{\rho^s-\overline{\rho^s}}_2
 \le C\norm{\nabla\rho^s}_2.
\]
It follows from \eqref{eq:basic} and \eqref{eq:bd} that
\begin{equation}\label{sy:Etime}
 \int_0^\infty \mathcal E(t)\dd t\le C,\qquad \lim_{t\to\infty}\mathcal E(t)=0.
\end{equation}
The limit follows from monotonicity.

For any fixed $0<d<R$, \eqref{eq:sphere} and the upper density
bound yield a uniform bound for $\rho$ in $H^1(d,R)$. Since
$\norm{\rho-\bar\rho}_{L^2(d,R)}\to0$ by \eqref{sy:Etime},
one-dimensional interpolation gives $\rho\to\bar\rho$ uniformly
on $[d,R]$. In particular, on this annulus $\rho\ge\bar\rho/2$
for all sufficiently large times. We use this lower bound only away
from the center in the next estimate.

\begin{lemma}\label{sy:angular-bound}
For every finite $q\ge2$, we have
\begin{equation}\label{sy:Wmoment}
 \sup_{t\ge0}\int\rho|W|^q\dd x\le C_q.
\end{equation}
Moreover, for every finite $m\ge2$,
\begin{equation}\label{sy:Ftime}
 \int_0^\infty\left(\int\rho F^m\dd x\right)^{1/m}dt\le C_m.
\end{equation}
\end{lemma}

\begin{proof}
Set $\omega=W/r$.
Multiply \eqref{sy:angular-stress} and \eqref{sy:angular} by
$|W|^{q-2}Wr^{2-q}$ and integrate with respect to $r\dd r$.
Adding the resulting identities with weights $1-1/q$ and $1/q$,
respectively, and using the continuity equation, we obtain
\begin{align*}
 \frac1q\frac{\dd}{\dd t}\int_0^R\rho r^3|\omega|^q\dd r
 &+\int_0^R\rho^\alpha r|\omega|^{q-2}
   \left[(q-1)(r\omega_r)^2+r\omega_r\omega+\frac2q\omega^2\right]\dd r\\
 &=-\int_0^R\rho r^2|\omega|^q
   \left[\left(2-\frac{5}{2q}\right)U+\frac{V}{2q}\right]\dd r .
\end{align*}
Here $W_r-W/r=r\omega_r$. 
For $q\ge2$, Young's inequality gives
\[
 (q-1)(r\omega_r)^2+r\omega_r\omega+\frac2q\omega^2
 \ge\frac{q-1}{2}(r\omega_r)^2+\frac1q\omega^2.
\]
Since $U,V\ge-A$ by Lemma \ref{sy:negative}, it follows that
\begin{equation}\label{sy:angular-diff}
 \frac1q\frac{\dd}{\dd t}\int_0^R\rho r^3|\omega|^q\dd r+c_q\int_0^R\rho^\alpha
  \bigl(r^3|\omega|^{q-2}\omega_r^2+r|\omega|^q\bigr)\dd r
 \le C_q\int_0^R\rho r^2|\omega|^q\dd r .
\end{equation}

By Lemma \ref{sy:upper}, $r\rho^{1-\alpha}\le Cr$. Choose $0<d<R$
so that $C_q\rho r^2\le(c_q/2)\rho^\alpha r$ on $[0,d]$.
On $[d,R]$ we have $\rho\ge\bar\rho/2$ for $t\ge t_d$ with $t_d$
sufficiently large, and hence $\int_d^R|\omega|^2\dd r\le C_d\mathcal E(t)$.
For $q>2$, applying the fundamental theorem of calculus to $|\omega|^{q/2}$
with $\omega(R,t)=0$ and then Young's inequality gives
\begin{align*}
 \int_d^R|\omega|^q\dd r
 &\le \norm{\omega}_{L^\infty(d,R)}^{q-2}\int_d^R|\omega|^2\dd r\\
 &\le C_{d,q}\left(\int_d^R|\omega|^{q-2}\omega_r^2\dd r\right)^{1-2/q}
          \mathcal E(t)\\
 &\le \varepsilon\int_d^R|\omega|^{q-2}\omega_r^2\dd r
       +C_{\varepsilon,d,q}\mathcal E(t)^{q/2}.
\end{align*}
For $q=2$, the $L^2(d,R)$ estimate suffices. Absorbing the derivative
term in \eqref{sy:angular-diff} for small $\varepsilon$ and using
$\rho r^3\le C\rho^\alpha r$, we obtain
\[
 \frac{\dd}{\dd t}\int_0^R\rho r^3|\omega|^q\dd r
 +c_q\int_0^R\rho r^3|\omega|^q\dd r
 \le C_q\mathcal E(t)^{q/2},\qquad t\ge t_d.
\]
Classical regularity bounds $W/r$ on $[0,2t_d]$. For $t\ge2t_d$,
integration from $t_d$ and the monotonicity of $\mathcal E$ give
\begin{align*}
 \int_0^R\rho r^3|\omega|^q\dd r
 &\le C_q e^{-c_q(t-t_d)}
   +C_q\int_{t_d}^t e^{-c_q(t-s)}\mathcal E(s)^{q/2}\dd s\\
 &\le C_q e^{-c_qt/2}+C_q\mathcal E(t/2)^{q/2},
\end{align*}
where we split the time integral at $t/2$. Since
$\int\rho|W|^q\dd x=2\pi\int_0^R\rho r^{q+1}|\omega|^q\dd r$
and $r^{q+1}\le R^{q-2}r^3$, this proves \eqref{sy:Wmoment}.

Finally, $F=r\omega^2$ and $r^{m+1}\le R^{m-2}r^3$ for $m\ge2$.
Taking $q=2m$ in the decay estimate gives, for all sufficiently large $t$,
\[
 \left(\int\rho F^m\dd x\right)^{1/m}
 \le \left(2\pi R^{m-2}\int_0^R\rho r^3|\omega|^{2m}\dd r\right)^{1/m}
 \le C_m\bigl(e^{-c_m t}+\mathcal E(t/2)\bigr).
\]
Integrating in time and using \eqref{sy:Etime} proves \eqref{sy:Ftime}.
\end{proof}

\subsection{The time-uniform lower bound of the density}

\begin{proposition}\label{sy:uniform}
For $1/2<\alpha<1$, let $2<m<N_\alpha$, where $N_\alpha$ is given
by \eqref{oa:N}. For $\alpha=1$, let $m>2$ be finite. Then
\begin{align}
 &\sup_{t\ge0}\int\rho\bigl(|u|^m+|V|^m\bigr)\dd x\notag\\
 &\quad+\int_0^\infty\int
 \left\{\rho^\alpha |U|^{m-2}
           \left(U_r^2+\frac{U^2}{r^2}\right)
          +\rho^{\gamma-\alpha+1}|V|^m\right\}\dd x\dd t\le C_m .
 \label{sy:radialmoments}
\end{align}
Moreover,
\begin{equation}\label{sy:uniformdensity}
 C^{-1}\le\rho(r,t)\le C
 \qquad (0\le r\le R,\ t\ge0),\qquad
 \lim_{t\to\infty}\norm{\rho(t)-\bar\rho}_\infty=0 .
\end{equation}
\end{proposition}

\begin{proof}
Multiplying the momentum equation in \eqref{eq:model} by $|U|^{m-2}Ue_r$
and using the continuity equation, we obtain
\begin{align*}
 &\frac1m\frac{\dd}{\dd t}\int\rho|U|^m\dd x
 +2\int\rho^\alpha|U|^{m-2}
 \left[\alpha(m-1)U_r^2-(1-\alpha)mU_r\frac Ur
              +\alpha\frac{U^2}{r^2}\right]\dd x\\
 &\qquad=-\int(\rho^\gamma)_r|U|^{m-2}U\dd x
           +\int\rho F|U|^{m-2}U\dd x .
\end{align*}
 For $2<m<N_\alpha$,
\eqref{oa:N} gives $4\alpha^2(m-1)>(1-\alpha)^2m^2$.
When $\alpha=1$, this holds for every $m>2$. Thus the quadratic form is bounded
below by $c_m(U_r^2+U^2/r^2)$.

To estimate the pressure, fix $K\ge1$ and split the multiplier as
\[
 |U|^{m-2}U
 =U\min\{|U|^{m-2},K^{m-2}\}+U(|U|^{m-2}-K^{m-2})_+ .
\]
For the first term, Cauchy--Schwarz and the upper density bound give
\begin{align*}
 \left|\int(\rho^\gamma)_r U\min\{|U|^{m-2},K^{m-2}\}\dd x\right|
 &\le \gamma K^{m-2}
   \left(\int\rho^{\gamma+\alpha-3}\rho_r^2\dd x\right)^{1/2}
   \left(\int\rho^{\gamma-\alpha+1}U^2\dd x\right)^{1/2}\\
 &\le C_{m,K}\left(\int\rho^{\gamma+\alpha-3}\rho_r^2\dd x
                         +\mathcal E(t)\right).
\end{align*}
Here $\rho^{\gamma-\alpha+1}\le C\rho$ since $\gamma>\alpha$.
For the second term, integrate by parts. The multiplier vanishes at
$|U|=K$ and at $r=R$. On $|U|>K$, we have $R|U|/(Kr)\ge1$. Hence
\begin{align*}
 \left|\int(\rho^\gamma)_r U(|U|^{m-2}-K^{m-2})_+\dd x\right|
 &\le C_m\int_{|U|>K}\rho^\gamma |U|^{m-2}
                 \left(|U_r|+\frac{|U|}{r}\right)\dd x\\*
 &\le \frac{C_m}{K}\int_{|U|>K}\rho^\alpha |U|^{m-2}
                 \left(\frac{|U|}{r}|U_r|+\frac{U^2}{r^2}\right)\dd x\\*
 &\le \frac{C_m}{K}\int\rho^\alpha |U|^{m-2}
                 \left(U_r^2+\frac{U^2}{r^2}\right)\dd x .
\end{align*}
Choose $K$ large enough to absorb this term into the viscous dissipation.
For the centrifugal term, H\"older's inequality gives
\[
 \left|\int\rho F|U|^{m-2}U\dd x\right|
 \le\left(\int\rho F^m\dd x\right)^{1/m}
    \left(\int\rho|U|^m\dd x\right)^{1-1/m}.
\]
Integrating in time and using \eqref{eq:bd}, \eqref{sy:Etime} and
\eqref{sy:Ftime}, we obtain for every $T>0$
\begin{equation}\label{sy:U-moment}
\begin{split}
 &\sup_{0\le t\le T}\int\rho|U|^m\dd x
   +\int_0^T\int\rho^\alpha |U|^{m-2}
                \left(U_r^2+\frac{U^2}{r^2}\right)\dd x\dd t\\
 &\qquad\le C_m+C_m\left(\sup_{0\le t\le T}\int\rho|U|^m\dd x\right)^{1-1/m}.
\end{split}
\end{equation}
By Young's inequality the last term is at most half the supremum on
the left plus $C_m$. This bounds both terms on the left independently of $T$.

Next, multiply \eqref{eq:radial-effective} by $|V|^{m-2}V$.
Using Young's inequality on the term containing $U$ and H\"older's
inequality on the term containing $F$, we find
\begin{align*}
 \frac1m\frac{\dd}{\dd t}\int\rho|V|^m\dd x
  +c\int\rho^{\gamma-\alpha+1}|V|^m\dd x
 &\le C_m\int\rho^{\gamma-\alpha+1}|U|^m\dd x\\
 &\quad+\left(\int\rho F^m\dd x\right)^{1/m}
         \left(\int\rho|V|^m\dd x\right)^{1-1/m}.
\end{align*}
Since $\gamma-2\alpha+1>0$, Lemma \ref{sy:upper} gives
\[
 \int\rho^{\gamma-\alpha+1}|U|^m\dd x
 \le R^2\norm{\rho}_\infty^{\gamma-2\alpha+1}
       \int\rho^\alpha\frac{|U|^m}{r^2}\dd x .
\]
The right-hand side is integrable in time by \eqref{sy:U-moment}.
For the centrifugal term, \eqref{sy:Ftime} and Young's inequality give
\[
 \int_0^T\left(\int\rho F^m\dd x\right)^{1/m}
         \left(\int\rho|V|^m\dd x\right)^{1-1/m}\dd t
 \le \frac1{2m}\sup_{0\le t\le T}\int\rho|V|^m\dd x+C_m .
\]
Integrating the $V$ inequality and absorbing the supremum gives both
$V$ estimates in \eqref{sy:radialmoments}. The estimate for $u$ follows
from \eqref{sy:Wmoment} and $|u|^m\le C_m(|U|^m+|W|^m)$.

Finally, \eqref{eq:effective} gives $2\alpha\rho_r=\rho^{2-\alpha}(V-U)$.
Since $m(2-\alpha)>1$ and $\rho$ is bounded above, we have
\[
 \int|\nabla\rho|^m\dd x
 \le C_m\int\rho^{m(2-\alpha)}\bigl(|V|^m+|U|^m\bigr)\dd x
 \le C_m\int\rho\bigl(|V|^m+|U|^m\bigr)\dd x\le C_m .
\]
By \eqref{sy:Etime} and \eqref{eq:relative-equivalence},
$\norm{\rho(t)-\bar\rho}_2\to0$. Since $m>2$, interpolation on the disk gives
\[
 \norm{\rho-\bar\rho}_\infty
 \le C_m\norm{\rho-\bar\rho}_2^{\frac{m-2}{2(m-1)}}
          \norm{\nabla\rho}_m^{\frac{m}{2(m-1)}}
      +C_m\norm{\rho-\bar\rho}_2\to0 .
\]
Hence $\rho\ge\bar\rho/2$ for $t\ge T_0$ with $T_0$
sufficiently large. On $[0,T_0]$, Proposition \ref{sy:finitetime}
bounds $\rho^{-1}$. This proves \eqref{sy:uniformdensity}.
\end{proof}

For later use, \eqref{sy:radialmoments} and
\eqref{sy:uniformdensity} yield
\begin{equation}\label{sy:inputs}
\begin{split}
 &\sup_{t\ge0}\bigl(\norm{u(t)}_m+
     \norm{\nabla\rho(t)}_2+\norm{\nabla\rho(t)}_m\bigr)\\
 &\quad+\int_0^\infty
   \bigl(\norm{\nabla\rho(t)}_2^2+\norm{\nabla\rho(t)}_m^m+
                     \norm{\nabla u(t)}_2^2\bigr)\dd t\le C_m .
\end{split}
\end{equation}
Indeed, the quadratic estimates follow from \eqref{eq:basic},
\eqref{eq:bd}, and \eqref{eq:sphere}, while
$\int_0^\infty\norm{\nabla\rho}_m^m\dd t$ is controlled by the
time integrals of $\rho^\alpha|U|^m/r^2$ and
$\rho^{\gamma-\alpha+1}|V|^m$ in \eqref{sy:radialmoments}.

\begin{proof}[Proof of Theorem \ref{thm:main}]
Lemma \ref{lem:local} and the continuation arguments above give a unique
global classical solution. The density bounds \eqref{eq:main-density}
follow from \eqref{oa:lowerestimate}, \eqref{ep:uniformdensity} and
\eqref{sy:uniformdensity} in the respective cases. With the bounds for
$m>2$ in \eqref{oa:bridgeestimate}, \eqref{ep:uniforminputs} and
\eqref{sy:inputs}, the higher-order estimates and decay argument of
\cite[Sections 3.6--3.8]{HMZ2025} yield \eqref{eq:main-asymptotic}.
\end{proof}

\section*{Acknowledgments}
Xiangdi Huang acknowledges partial support from the Chinese Academy of
Sciences Project for Young Scientists in Basic Research (Grant No.~YSBR-031),
the National Natural Science Foundation of China (Grant Nos.~12494542 and
11688101), and the National Key R\&D Program of China (Grant
No.~2021YFA1000801).

\medskip
\noindent\textbf{AI assistance.}
Generative AI contributed substantially to the construction of the most important
mathematical proofs in this paper. Two particularly insightful examples among them
are the coefficient-freezing argument for higher integrability of the radial
effective velocity (gradient stress, $\alpha<1$, Proposition~\ref{oa:lower})
and the time-uniform, density-weighted $L^p$ bounds for the tangential velocity
(symmetric stress, Lemma~\ref{sy:angular-bound}). The authors also used AI to
explore proof strategies, check calculations, and polish the manuscript.
They take full responsibility for the mathematics and the final text.

\medskip
\noindent\textbf{Data availability statement.}
No data were generated or analyzed in this study.

\medskip
\noindent\textbf{Conflict of interest.}
The authors declare no competing interests.


\end{document}